\input amstex
\magnification=\magstep1 
\baselineskip=13pt
\documentstyle{amsppt}
\vsize=8.7truein \CenteredTagsOnSplits \NoRunningHeads
\def\today{\ifcase\month\or
  January\or February\or March\or April\or May\or June\or
  July\or August\or September\or October\or November\or December\fi
  \space\number\day, \number\year}
 \def\per{\operatorname{per}}
 \def\EE{\bold{E}\thinspace}
 \def\PP{\bold{P}\thinspace}

 \def\PPP{\Cal P}
 \topmatter
 \abstract
 We present a randomized approximation algorithm for counting {\it contingency tables},
 $m \times n$ non-negative 
 integer matrices with given row sums $R=\left(r_1, \ldots, r_m\right)$
 and column sums $C=\left(c_1, \ldots, c_n\right)$.  We define
 {\it smooth margins} $(R,C)$ in terms of the {\it typical table} 
 and prove that for such margins the algorithm has quasi-polynomial 
$N^{O(\ln N)}$ complexity, where $N=r_1+\cdots+r_m=c_1+\cdots+c_n$.
 Various classes of margins are smooth, e.g., when $m=O(n)$, $n=O(m)$ and 
 the ratios between the largest and the smallest row sums as well as  
 between the largest and the smallest column sums are strictly 
 smaller than the golden ratio $(1+\sqrt{5})/2 \approx 1.618$. The algorithm builds 
 on Monte Carlo integration and sampling algorithms for log-concave densities, 
the matrix scaling algorithm, the permanent approximation algorithm, and an integral representation for the number of contingency tables.
 \endabstract
\title An approximation algorithm for counting contingency tables \endtitle
\author  Alexander Barvinok, Zur Luria,  Alex Samorodnitsky, and  Alexander Yong\endauthor
\address Department of Mathematics, University of Michigan, Ann Arbor,
MI 48109-1043, USA \endaddress
\email barvinok$\@$umich.edu \endemail

\address Department of Computer Science, Hebrew University of Jerusalem, Givat Ram Campus,
91904, Israel \endaddress
\email zluria$\@$cs.huji.ac.il \endemail
\address Department of Computer Science, Hebrew University of Jerusalem, Givat Ram Campus,
91904, Israel \endaddress
\email salex$\@$cs.huji.ac.il \endemail
\address Department of Mathematics, University of Minnesota, Minneapolis, MN 55455, USA
\endaddress
\email ayong$\@$math.umn.edu \endemail
\date March 2008 \enddate
\endtopmatter
\document

\head 1. Introduction  \endhead 

Let $R=\left(r_1, \ldots, r_m\right)$ and $C=\left(c_1, \ldots, c_n \right)$ be positive 
integer vectors such that 
$$\sum_{i=1}^m r_i =\sum_{j=1}^n c_j =N.$$
A {\it contingency table} with {\it margins} $(R, C)$ is an $m \times n$ non-negative 
integer matrix $D=\left(d_{ij}\right)$ with row sums $R$ and column sums $C$:
$$\split &\sum_{j=1}^n d_{ij}=r_i \quad \text{for} \quad i=1, \ldots, m \\
&\sum_{i=1}^m d_{ij} =c_j \quad \text{for} \quad j=1, \ldots, n. \endsplit$$
Let $\#(R,C)$ denote the number of these contingency tables.

There is interest in the study of $\#(R,C)$, due to 
connections to statistics, combinatorics and representation theory, see, e.g.,  
\cite{Go76},
\cite{DE85}, \cite{DG95},
\cite{D+97}, \cite{Mo02}, \cite{CD03}, \cite{L+04}, \cite{B+04}, \cite{C+05} and the
references therein. However, since enumerating $\#(R,C)$ is 
a $\#P$-complete problem even for $m=2$ \cite{D+97}, one does
not expect to find polynomial-time algorithms (nor formulas)
computing $\#(R,C)$ exactly. As a result, attention has turned to the 
open problem of efficiently estimating $\#(R,C)$. 

We present a randomized algorithm for approximating $\#(R,C)$ 
within a prescribed relative error.
Based on earlier numerical studies \cite{Yo07} \cite{B+07}, we
conjecture that its complexity is polynomial in $N$.
We provide further evidence for this hypothesis: we
introduce ``smooth margins'' $(R,C)$ where
the entries of the {\it typical table} are not too large, and 
among $\{r_1,\ldots,r_m,c_1,\ldots,c_n\}$ there are 
no ``outliers''. Our main result is that 
smoothness implies a quasi-polynomial 
$N^{O(\log N)}$ complexity bound on the algorithm. 
More precisely, we approximate $\#(R,C)$ within relative error $\epsilon>0$
using $(1/\epsilon)^{O(1)}N^{O(\ln N)}$ time in the unit cost model, provided 
$\epsilon \gg 2^{-m}+2^{-n}$.\footnote{If an exponentially small relative error $\epsilon=O\left(2^{-m} +2^{-n}\right)$ is desired, one has an exact dynamic programming algorithm with
$N^{O(m+n)}=(1/\epsilon)^{O(\ln N)} $ quasi-polynomial complexity.} 

The class of smooth margins captures a number of interesting subclasses. 
In particular, this work applies to the case of {\it magic squares} (where
$m=n$ and $r_i=c_j=t$ for all $i,j$), extending \cite{B+07}. 
More generally, smoothness includes the case when the ratios $m/n$ and $n/m$ are 
bounded by a constant fixed in advance while the ratios between the largest  
and the smallest row sums as well as between the largest and the smallest column sums 
are smaller than the golden ratio $\left(1+\sqrt{5}\right)/2 \approx 1.618$. 
These and others examples are explicated in Section~3. See Section~1.4 for 
comparisons to the literature.

\subhead (1.1) An outline of the algorithm \endsubhead 
Our algorithm builds on the 
technique of rapidly mixing Markov chains and, in particular, on efficient 
integration and sampling from log-concave densities, 
as developed in \cite{AK91}, \cite{F+94}, \cite{FK99}, \cite{LV06} (see also \cite{Ve05} for a survey), the permanent approximation algorithm
\cite{J+04}, the strongly polynomial time algorithm for matrix scaling \cite{L+00}, and 
the integral representation of $\#(R,C)$ from \cite{Ba08}.

Let $\Delta=\Delta_{m \times n} \subset {\Bbb R}^{mn}$ be the open $(mn-1)$-dimensional simplex of all $m \times n$ positive matrices $X=\left(x_{ij}\right)$
such that 
$$\sum_{ij} x_{ij}=1.$$
 Let $dX$ be Lebesgue measure
on $\Delta$ normalized to the probability measure. An integral
representation for $\#(R,C)$ was found in \cite{Ba08}: 
$$\#(R,C)=\int_{\Delta} f(X) \ d X, \tag1.1.1$$
where $f: \Delta \longrightarrow {\Bbb R}_+$ is a certain
continuous function that factors as 
$$f=p \phi, \tag1.1.2$$
where 
$$p(X) \geq 1 \quad \text{for all} \quad X \in \Delta$$
is a function that ``does not vary much'',
and
$\phi: \Delta \longrightarrow {\Bbb R}_+$ is  continuous 
and {\it log-concave}, that is,
$$\split \phi(\alpha X +\beta Y) \geq \phi^{\alpha}(X) \phi^{\beta}(Y) \quad& \text{for all} \quad 
X, Y \in \Delta  \quad \text{and} \\ & \text{for all} \quad \alpha,\beta \geq 0 \quad \text{such that} \quad 
\alpha+\beta=1. \endsplit$$
Full details about $f$ and its factorization are reviewed in Section~2.

For any $X \in \Delta$, the values of $p(X)$ and $\phi(X)$ are computable 
in time polynomial
in $N$. 
Given $\epsilon>0$,
the value of $p(X)$ can be computed, 
within relative error $\epsilon$ in time polynomial
in $1/\epsilon$ and $N$, by a randomized algorithm of \cite{J+04}.
The value of $\phi(X)$ can be 
computed, within relative error $\epsilon$ in time polynomial in $\ln (1/\epsilon)$ 
and $N$, by a deterministic 
algorithm of \cite{L+00}. 

The central idea of this paper is to define smooth margins $(R,C)$ so that
matrices $X \in \Delta$ with large values of $p(X)$ do not contribute much to the 
integral (1.1.1). Our main results, precisely stated in Section~3, are
that for smooth margins, there is a {\it threshold} $\tau=N^{\delta \ln N}$ for some 
constant $\delta>0$ (depending on the class of margins considered) such that 
if we define the truncation $\overline{p}: \Delta \longrightarrow {\Bbb R}_+$ by 
$$\overline{p}(X)=\cases p(X) &\text{if \ } p(X) \leq \tau \\ \tau &\text{if \ } p(X) >\tau \endcases$$
then 
$$\#(R, C)= \int_{\Delta} p(X) \phi(X) \ d X  \approx \int_{\Delta} \overline{p}(X) \phi(X) \ d X \tag1.1.3$$
where ``$\approx$'' means ``approximates
to within an $O\left(2^{-n}+2^{-m}\right)$ relative error'' 
(in fact, rather than base $2$, any constant $M>1$, fixed in advance,
can be used). 
We conjecture that one can choose the threshold $\tau =N^{O(1)}$, which would make the complexity of our algorithm polynomial in $N$.

The first step (and a simplified version) of our algorithm computes 
the integral $$\int_{\Delta} \phi(X) \ d X \tag1.1.4$$
using any of the aformentioned randomized polynomial time algorithms 
for {\it integrating} log-concave densities; these results imply that
this step has polynomial in $N$ 
complexity. By (1.1.3) it follows that for smooth $(R,C)$ the integral
(1.1.4) approximates $\#(R,C)$ within a factor of $N^{O(\ln N)}$. This simplified
algorithm is suggested in \cite{Ba08}; an implementation that utilizes a version of 
the hit-and-run algorithm of \cite{LV06}, 
together with numerical results is described in \cite{Yo07} and \cite{B+07}.

Next, our algorithm estimates (1.1.3) within relative error 
$\epsilon$ using the aformentioned randomized polynomial time algorithm 
for approximating the permanent of a matrix, and any of those for
{\it sampling} from log-concave densities. Specifically,
let $\nu$ be the probability measure on $\Delta$ with the 
density proportional to~$\phi$. Thus, 
$$\int_{\Delta} \overline{p}(X) \phi(X) \ dX=
\left(\int_{\Delta} \overline{p} \ d \nu  \right) \left( \int_{\Delta} \phi(X) \  d X \right).$$
The second factor is computed by the above first step, while
the first factor is approximated by the sample mean
$$\int_{\Delta} \overline{p} \ d \nu \approx {1 \over k} \sum_{i=1}^k \overline{p}(X_i), \tag1.1.5$$
where $X_1, \ldots, X_k \in \Delta$ are independent points sampled at random from measure $\nu$. 
Since $1 \leq \overline{p}(X) \leq \tau$, 
the Chebyshev inequality implies that to achieve relative error $\epsilon$
with probability $2/3$ it suffices to sample $k=O\left(\epsilon^{-2} \tau^2 \right)
=\epsilon^{-2} N^{O(\ln N)}$
points in (1.1.5).

The results of \cite{AK91}, \cite{F+94}, \cite{FK99}, 
and \cite{LV06} imply that 
for any given $\epsilon>0$ one can sample independent points $X_1, \ldots, X_k$ from  
a distribution $\tilde{\nu}$ on $\Delta$ such that 
$$|\tilde{\nu}(S)-\nu(S)| \leq \epsilon \quad \text{for any Borel set} \quad S \subset \Delta.$$
in time linear in $k$ and polynomial in $\epsilon^{-1}$ and $N$.
Replacing $\nu$ by $\tilde{\nu}$ in (1.1.5) introduces an additional relative error 
of $\epsilon \tau =\epsilon N^{\delta \ln N}$, handled by choosing a smaller 
$\epsilon =O\left(N^{-\delta \ln N}\right)$.

\subhead (1.2) An optimization problem, typical tables and smooth margins\endsubhead 
We will define smoothness of margins in terms of a certain convex optimization problem.

Let $\PPP=\PPP(R,C)$ be the {\it transportation polytope} of $m \times n$ non-negative
matrices $X=\left(x_{ij}\right)$ with row sums $R$ and column sums $C$.
On the space ${\Bbb R}^{mn}_+$ of $m \times n$ 
non-negative matrices define
$$g(X)=\sum_{ij} \Bigl( \left(x_{ij}+1\right) \ln \left(x_{ij}+1\right) -x_{ij} \ln x_{ij} \Bigr)
\quad \text{for} \quad X=\left(x_{ij}\right).$$
The following optimization problem plays an important role in this paper:
$$ \text{Maximize} \quad g(X) \quad \text{subject to} \quad X\in \PPP.\tag1.2.1$$
It is easy to check that $g$ is strictly concave and hence attains its maximum on $\PPP$
at a unique matrix $X^{\ast}=\left(x_{ij}^{\ast}\right)$, $X^{\ast} \in \PPP$ that we call
the {\it typical table}.

An intuitive explanation for the appearance of this optimization problem,
and justification for the nomenclature ``typical'' derives from work of 
\cite{B07b} (relevant parts are replicated for convenience, 
in Section 4, see specifically Theorem 4.1). In short, 
$X^{\ast}$ determines the asymptotic behavior of $\#(R,C)$.

The main requirement that we demand of 
smooth margins $(R, C)$ to satisfy (see Section~3 
for unsuppressed technicalities) is that the entries of 
the typical table are not too large, that is, entries $x_{ij}^{\ast}$ of the 
optimal solution $X^{\ast}=\left(x_{ij}^{\ast}\right)$ satisfy
$$\max_{ij} x_{ij}^{\ast} = O(s) \quad \text{where} \quad s={N \over mn}$$
is the average entry of the table.

Viewing the typical table as interesting in its own right, 
one would like to understand 
how the typical table changes as the margins vary.
The optimization problem being convex, $X^{\ast}$ can be computed efficiently 
by many existing algorithms, see, for example, \cite{NN94}. However, 
in many instances of interest, the smoothness
condition can be checked without actually needing to solve this problem.
For example, if all the row sums 
$r_i$ are equal, the symmetry of the functional $g$ under permutations of rows implies that
$$x_{ij}^{\ast}={c_j \over m} \quad \text{for all} \quad i,j.$$
In general, the entries $x_{ij}^{\ast}$ stay small if the row sums $r_i$ and 
column sums $c_j$ do not vary much. On the other hand,
it is not hard to construct examples of margins $(R, C)$ for $n$-vectors 
$R$ and $C$ such that $n \leq r_i, c_j \leq 3n$ and some of the entries $x_{ij}^{\ast}$ are large, 
in fact linear in $n$.  Another one of our results (Theorem~3.5) gives upper and lower
bounds for $x_{ij}^{\ast}$ in terms of $(R,C)$. 

\subhead (1.4) Comparisons with the literature \endsubhead
Using the Markov Chain Monte Carlo approach, Dyer, Kannan, and Mount \cite{D+97} 
count contingency tables when $R$ and $C$ are sufficiently
large, that is, 
if $r_i=\Omega\left(n^2 m\right)$ and $c_j=\Omega\left(m^2 n \right)$ for all $i,j$. 
Their randomized (sampling) algorithm approximates $\#(R,C)$ 
within any given relative error $\epsilon>0$ 
in time polynomial in $\epsilon^{-1}$, $n$, $m$, and 
$\sum_i \log r_i + \sum_j \log c_j$ (the 
bit size of the margins). Subsequently, Morris \cite{Mo02}
obtained a similar result for the bounds 
$r_i=\Omega\left(n^{3/2} m \ln m \right)$ and $c_j=\Omega\left(m^{3/2} n \ln n\right)$. 
These results are based on fact that
for large margins, the number of contingency tables is well-approximated by the
volume of the transportation polytope $\PPP(R,C)$ (contingency tables being the
integer points in this polytope). More generally, Kannan and Vempala 
\cite{KV99} show
that estimating the number integer points in a $d$-dimensional 
polytope with $m$ facets reduces to 
computing the volume of the polytope (a problem, for which efficient 
randomized algorithms
exist, see \cite{Ve05} for a survey) provided the polytope contains a ball of radius $d \sqrt{\log m}$.

When the margins $r_i, c_j$ are very small, that is, 
bounded by a constant fixed in advance) relative to the 
sizes $m$ and $n$ of the matrix,  B\' ek\' essy, 
B\' ek\' essy, and Koml\' os \cite{B+72} obtain an 
efficient and precise asymptotic formula for $\#(R,C)$.
Their formula exploits the fact in this case, the majority of contingency tables 
have only entries $0$, $1$, and $2$. Alternatively, in this case one can 
exactly compute $\#(R,C)$ in time polynomial in $m+n$ via a dynamic programming 
algorithm. More recently, Greenhill and McKay \cite{GM07} gave 
a computationally efficient 
asymptotic formula for a wider class of sparse 
margins (when $r_i c_j =o(N^{2/3}) )$. 

Also using the dynamic programming approach, 
Cryan and Dyer \cite{CD03} construct a randomized polynomial time approximation
algorithm to compute $\#(R,C)$, provided the number of rows is fixed; 
see \cite{C+06} for sharpening of the results. 

It seems that the most resilient case of computing $\#(R,C)$ is
where both $m$ and $n$ grow, and the margins are of moderate size, e.g., 
linear in the dimension. 
Recently, Canfield and McKay \cite{CM07} 
found a precise asymptotic formula for $\#(R,C)$ assuming that all row sums are 
equal and all column sums are equal. However, for general margins no such 
formula is known, even conjecturally.

We remark that our notion of smooth margins includes all of the above 
regimes, except for that of large margins.

Summarizing, although our  complexity bounds do not improve on the
algorithms in the above cases, our algorithm is provably computationally efficient 
(quasi-polynomial in $N$) for several new classes of margins, which include cases 
of growing dimensions $m$ and $n$ and moderate size margins $R$ and $C$. 

\head 2. The integral representation for the number of contingency tables \endhead

We now give details of the integral representation (1.1.1). To do this, 
we express $\#(R,C)$ as the 
expectation of the permanent of a random $N \times N$  matrix. Recall that the 
{\it permanent} of an $N \times N$ matrix $A$ is defined by 
$$\per A=\sum_{\sigma  \in S_N} \prod_{i=1}^N a_{i \sigma(i)},$$
where $S_N$ is the symmetric group of the permutations of the set $\{1, \ldots, N\}$.
The following result was proved in \cite{Ba08}.
\proclaim{(2.1) Theorem} 
For an $m \times n$ matrix $X=\left(x_{ij}\right)$, let $A(X)$ be 
the $N \times N$ block matrix $A(X)$ whose the $(i,j)$-th block is the
$r_i \times c_j$ submatrix filled with $x_{ij}$, for $i=1, \ldots, m$ and
$j=1, \ldots, n$. 
Then
$${\per A(X) \over r_1! \cdots r_m! c_1! \cdots c_n!}
=\sum \Sb D=\left(d_{ij}\right) \endSb \prod_{ij} {x_{ij}^{d_{ij}} \over d_{ij}!},\tag2.1.1$$
where the sum is over all non-negative integer matrices $D=\left(d_{ij}\right)$
with row sums $R$ and column sums $C$.

Let ${\Bbb R}^{mn}_+$ be the open orthant of positive $m \times n$ matrices $X$.
Then
$$\#(R,C) ={1 \over r_1! \cdots r_m! c_1! \cdots c_n!}
\int_{{\Bbb R}^{mn}_+} \per A(X) \ \exp\left\{-\sum_{ij} x_{ij} \right\} \ dX,$$
where $dX$ is the Lebesgue measure on ${\Bbb R}^{mn}_+$.
\endproclaim

In the case that $r_i=a$ and $c_j=b$ for all $i,j$, the expansion (2.1.1) was first 
observed by Bang and then used by Friedland  \cite{Fr79} in his proof of a weaker 
form of the van der Waerden conjecture; see Section 7.1 and references there.

Since the function $X \longmapsto \per A(X)$ is a homogeneous polynomial of degree $N$, one can express
$\#(R, C)$ as an integral over the simplex. The following corollary was also obtained in 
\cite{Ba08}.

\proclaim{(2.2) Corollary} Let $\Delta=\Delta_{m \times n} \subset {\Bbb R}^{mn}$ be the open simplex of positive $m \times n$
matrices $X=\left(x_{ij}\right)$ such that $\sum_{ij} x_{ij}=1$. Then 
$$\#(R, C)={(N+mn-1)! \over (mn-1)!} {1 \over r_1! \ldots r_m! c_1! \cdots c_n!} \int_{\Delta_{m \times n}}
\per A(X) \ dX,$$
where $dX$ is the Lebesgue measure on $\Delta_{m \times n}$ normalized to the probability 
measure.
\endproclaim

Hence in the integral representation (1.1.1), 
we define the function $f$ by
$$\split f(X)=&{(N+mn-1)! \over (mn-1)!} {1 \over r_1! \ldots r_m! c_1! \cdots c_n!} \per A(X)\\
                      =&{(N+mn-1)! \over (mn-1)!} \sum \Sb D=\left(d_{ij}\right) \endSb \prod_{ij} 
                       {x_{ij}^{d_{ij}} \over d_{ij}!} , \endsplit$$
where $A(X)$ is the block matrix of Theorem 2.1 and the sum is 
over all contingency tables $D$ with margins $(R, C)$.

\subhead (2.3) Matrix scaling and the factorization of $f$ \endsubhead To 
obtain the factorization (1.1.2), where $\phi: \Delta \longrightarrow {\Bbb R}_+$ is 
a log-concave function and $p: \Delta \longrightarrow {\Bbb R}_+$ is a function which 
``does not vary much'', we employ the idea of {\it matrix scaling}, see
\cite{Si64}, \cite{MO68}, \cite{KK96}, Chapter 6 of  \cite{BR97}, and \cite{L+00}: 
Let $X=\left(x_{ij}\right)$ be a positive $m \times n$ matrix. Then there exists a unique 
$m \times n$ matrix $Y$ with the row sums $R=\left(r_1, \ldots, r_m \right)$, column
sums $C=\left(c_1, \ldots, c_n\right)$, and such that 
$$x_{ij}=y_{ij} \lambda_i \mu_j \quad \text{for all} \quad i,j$$
and some positive $\lambda_1, \ldots, \lambda_m, \mu_1, \ldots, \mu_n$.
The numbers $\lambda_i$ and $\mu_j$ are unique up to a re-scaling 
$ \lambda_i \longmapsto \lambda_i \tau$, $\mu_j \longmapsto \mu_j \tau^{-1}$.
Note that if we divide the entries in the $(i,j)$-th block of the matrix $A(X)$ of Theorem 2.1 by 
$r_i c_j \lambda_i \mu_j$, we obtain a positive {\it doubly stochastic} matrix $B(X)$, that is,
a positive matrix with all row and column sums equal to 1.
Thus we have 
$$\per A(X) = \left(\prod_{i=1}^m (\lambda_i r_i)^ {r_i} \right) 
\left(\prod_{j=1}^n (\mu_j c_j)^{c_j} \right)\per B(X) .$$
It is proved in \cite{Ba08} that 
$${N! \over N^N} \leq \per B(X) \leq \min\left\{ \prod_{i=1}^m {r_i! \over r_i^{r_i}}, \quad 
\prod_{j=1}^n {c_j! \over c_j^{c_j}} \right\}.$$
The lower bound is the van der Waerden bound for permanents of doubly stochastic 
matrices, see \cite{Fa81}, \cite{Eg81} and also Chapter 12 of \cite{LW01} and recent 
\cite{G06a}, while the upper bound is a corollary of the Minc conjecture proved 
by Bregman, see \cite{Br73},  Chapter 11 of \cite{LW01}, and also \cite{So03}.

Now we define 
$$p(X) ={N^N \over N!} \per B(X) \tag2.3.1$$ 
and 
$$\phi(X) ={(N+mn-1)! N! \over (mn-1)! N^N}  
\left(\prod_{i=1}^m {r_i^{r_i} \over r_i!} \right) \left(\prod_{j=1}^n {c_j^{c_j} \over c_j!} \right)
\left(\prod_{i=1}^m \lambda_i^{r_i} \right) \left(\prod_{j=1}^n \mu_j^{c_j} \right).$$
We summarize results of \cite{Ba08} regarding $p$ and $\phi$.

\proclaim{(2.4) Theorem} The following hold: 
\roster 
\item $\phi$ is log-concave, that is, 
$$\phi(\alpha X + \beta Y) \geq \phi^{\alpha}(X) \phi^{\beta}(Y)$$
for all $X, Y \in \Delta$ and $\alpha, \beta \geq 0$ such that $\alpha+\beta=1$;
\item Let $X, Y \in \Delta$ be positive $m \times n$ matrices, $X=\left(x_{ij}\right)$ and 
$Y=\left(y_{ij}\right)$, such that $x_{ij}, y_{ij} \geq \delta$ for all $i,j$ and some $\delta>0$.
Then 
$$\big|\ln \phi(X)-\ln \phi(Y)\big| \leq {N \over \delta} \max_{ij} \big| x_{ij}-y_{ij}\big|;$$
\item For $\delta<1/mn$ let us define the $\delta$-interior $\Delta_{\delta}$ of the 
simplex $\Delta$ as the set of matrices $X \in \Delta$, $X=\left(x_{ij}\right)$, such that 
$x_{ij} \geq \delta$ for all $i,j$. Then for  $f=p \phi$ we have 
$$(1-mn \delta )^{N+mn-1} \int_{\Delta} f \ dX \ \leq \ \int_{\Delta_{\delta}} f \ d X \ \leq \ 
\int_{\Delta} f \ dX;$$
\item We have 
$$1 \ \leq \ p(X) \ \leq \  {N^N \over N!} \min \left\{ \prod_{i=1}^m {r_i! \over r_i^{r_i}}, \ 
\prod_{j=1}^n {c_j! \over c_j^{c_j}} \right\}.$$
\endroster
\endproclaim
The log-concavity of function $\phi$ was first observed in \cite{G06b}.
In terms of \cite{G06b}, up to a normalization factor, $\phi(X)$ is the {\it capacity} 
of the matrix $A(X)$ of Theorem 2.1, see also \cite{B07b} for a more general
family of inequalities satisfied by $\phi$. As is discussed in \cite{Ba08}, the matrix 
scaling algorithm of \cite{L+00} leads to a polynomial time algorithm for computing 
$\phi(X)$. Namely, for any given $\epsilon>0$ the value of $\phi(X)$ can be 
computed within relative error of $\epsilon$ in time polynomial in $N$ and $\ln (1/\epsilon)$
in the unit cost model; our own experience is that this algorithm for computing
$\phi(X)$ is practical, and 
works well for $m,n \leq 100$. 

Theorems 2.4 and 2.1 allow us to apply algorithms of 
\cite{AK91}, \cite{F+94}, \cite{FK99},
and \cite{LV06} on efficient integration and sampling of log-concave functions. 
First, for any given $\epsilon>1$, one can compute the integral 
$$\int_{\Delta} \phi \ dX$$
within relative error $\epsilon$ 
in time polynomial in $\epsilon^{-1}$ and $N$ by a randomized algorithm. 
Second, one can sample points $X_1, \ldots, X_k \in \Delta$ independently from
a measure $\tilde{\nu}$ such that  
$$|\nu(S) - \tilde{\nu}(S)| \leq \epsilon \quad \text{for any Borel set} \quad S \subset \Delta,$$
where $\nu$ is the measure with the density proportional to $\phi$, in time polynomial 
in $k, \epsilon^{-1}$ and $N$.

The integration of $p(X)$ raises a greater challenge. 
For any given $\epsilon>0$ one can 
compute $p(X)$ itself 
within relative error $\epsilon$ in time polynomial in $\epsilon^{-1}$ and 
$N$, using the permanent approximation algorithm of \cite{J+04}. However,
the upper bound of Part (4) of Theorem 2.4 is, in the worst case, of order 
$N^{\gamma(m+n)}$ for some absolute constant $\gamma>0$. 
Therefore, {\it a priori}, to integrate $p$
over $\Delta$ using a sample mean, one needs too many such computations to guarantee
the desired accuracy of $\epsilon$. Our main observation to overcome this problem
is that in many interesting cases the matrices 
$X \in \Delta$ with large values of 
$p(X)$ do not contribute much to the integral (1.1.1), so we have $p(X)=N^{O(\ln N)}$
with high probability with respect to the density on $\Delta$ proportional to $f$.

\subhead (2.5) Bounding $p$ with high probability \endsubhead 
Let us consider the projection 
$$\split pr: {\Bbb R}^{mn}_+ \longrightarrow \Delta_{m \times n}, \quad pr(X)=&\tilde{X}, \quad \text{where}\\
\tilde{X} =\alpha X \quad \text{for} \quad \alpha =\left(\sum_{ij} x_{ij} \right)^{-1}. \endsplit$$
Clearly, the scalings of $X$ and $\tilde{X}$ to the matrix with the row sums $R$ and column
sums $C$ coincide. Also, it is clear that the doubly stochastic scalings 
$B(X)$ and $B(\tilde{X})$,   
of matrices $A(X)$ and $A(\tilde{X})$, respectively, also coincide. We define 
$p(X)$ for an arbitrary positive $m \times n$ matrix $X$ by 
$p(X):=p(\tilde{X})$, or, equivalently, by (2.3.1).

We introduce the following density $\psi=\psi_{R,C}$ on ${\Bbb R}^{mn}_+$ by
$$\split \psi(X)={1 \over \#(R, C)} \sum \Sb D=\left(d_{ij}\right) \endSb \prod_{ij} &{x_{ij}^{d_{ij}} \over d_{ij}!}
e^{-x_{ij}}, \quad \text{where} \\ &X=\left(x_{ij}\right) \quad \text{and} \quad x_{ij} >0
\quad \text{for all} \quad i,j, \endsplit$$
and the sum is over all $m \times n$ non-negative integer matrices $D$ with the 
row sums $R$ and column sums $C$.
We define $\psi(X)=0$ if $X$ is not a positive matrix. 
That $\psi$ is a probability density is immediate from Theorem~2.1. 

Our goal 
is to show that for smooth margins $(R,C)$, the value of 
$p(X)$ is ``reasonably small'' for most $X$, that is,
 $$\PP\left\{X \in {\Bbb R}^{mn}_+: \quad p(X) > N^{\delta \ln N} \right\}  < \kappa \left(2^{-m}+2^{-n}\right) \tag2.5.1$$
for some constants $\delta>0$ and $\kappa>0$, where the probability is measured with respect to the density $\psi$.

Our construction of function $f$ in (1.1.1) implies that the push-forward of $\psi$ under the projection
$pr: {\Bbb R}^{mn}_+ \longrightarrow \Delta$ 
is the density 
$${1 \over \#(R, C)} f(X) \quad \text{for} \quad X \in \Delta$$
on the simplex. Hence inequality (2.5.1) implies that for $\tau=N^{\delta \ln N}$ we have 
$${1 \over \#(R, C)} \int\limits \Sb X \in \Delta \\  p(X) > \tau \endSb  f(X) \ dX < \kappa \left(2^{-m} + 2^{-n} \right).$$
Therefore, as discussed in Section 1.1, replacing $p$ by its truncation $\overline{p}$ introduces an 
$O\left(2^{-n} + 2^{-m}\right)$ relative error in (1.1.3) and
hence our algorithm achieves quasi-polynomial complexity.

The key idea behind inequality (2.5.1) is that the permanent of an 
appropriately defined ``random'' doubly stochastic matrix is very close with 
high probability to the van der Waerden lower bound $N!/N^N$; see Lemma~5.1.

\head 3. Main results \endhead 

Now we are ready to precisely define the classes of smooth margins
for which our algorithm achieves $N^{O(\ln N)}$ complexity. 

\definition{ (3.1) Smoothness Definitions}
Fix margins $R=\left(r_1, \ldots, r_m\right)$, $C=\left(c_1, \ldots, c_n \right)$,
where 
$$\sum_{i=1}^m r_i =\sum_{j=1}^n c_j =N.$$
Let
$$s={N \over mn}$$
be the average value of the entries of the table.
We define
$$\split &r_+=\max_{i=1, \ldots, m} r_i, \quad r_-=\min_{i=1, \ldots, m} r_i \\
&c_+=\max_{j=1, \ldots, n} c_j, \quad c_-=\min_{j=1, \ldots, n} c_j. \endsplit$$
Hence $r_+$ and $c_+$ are the largest row and column sums respectively and
$r_-$ and $c_-$ are the smallest row and column sums respectively.

For $s_0>0$, call the margins $(R, C)$ {\it $s_0$-moderate} if 
$s \leq s_0$.
In other words, margins are moderate if the average entry of the table is bounded from 
above.

For $\alpha \geq 1$, the margins $(R, C)$ are {\it upper $\alpha$-smooth} if
$$r_+ \ \leq \ \alpha sn=\alpha {N \over m}  \quad \text{and} \quad c_+\  \leq \ \alpha sm=\alpha 
{N \over n}.$$
Thus, margins are upper smooth if the row and column sums are at most proportional 
to the average row and column sums respectively.

For $0< \beta \leq 1$, the margins $(R,C)$ are {\it lower $\beta$-smooth}
if 
$$r_-\  \geq\ \beta sn=\beta {N \over m} \quad \text{and} \quad c_- \ \geq \ \beta sm=\beta {N \over n}.$$
Therefore, margins are lower smooth if the row and column sums are at least 
proportional to the average row and column sums respectively.

The key smoothness condition is as follows: 
for $\alpha \geq 1$ we define margins $(R,C)$ to be 
{\it strongly upper $\alpha$-smooth}
if for the typical table $X^{\ast}=\left(x_{ij}^{\ast}\right)$ we have
$$x_{ij}^{\ast} \leq \alpha s \quad \text{for all} \quad i,j.$$
Note that this latter condition implies that the margins are upper $\alpha$-smooth.
(Also, we do not need a notion of strongly lower $\beta$ smooth.)
\enddefinition

Our main results are randomized approximation algorithms of quasi-polynomial 
$N^{O(\ln N)}$ complexity when the margins $(R,C)$ are smooth for either:
\medskip
$\bullet$ $s_0$-moderate strongly upper $\alpha$-smooth, for some fixed $s_0$ and $\alpha$;
\smallskip
\noindent or
\smallskip
$\bullet$ lower $\beta$ and strongly upper $\alpha$-smooth, for some fixed $\alpha$ and $\beta$.
\medskip

By the discussion of Section~2.5, the quasi-polynomial complexity claim about our
algorithm follows from bounding on $p(X)$ with high probability. Specifically, we 
have the following two results. Their proofs are argued similarly, but the second
is more technically involved.

\proclaim{(3.2) Theorem} Fix $s_0 >0$ and $\alpha \geq 1$. Suppose that $m \leq 2^n$, 
$n \leq 2^m$ and
let $(R,C)$ be $s_0$-moderate strongly upper $\alpha$-smooth margins.
Let $X=\left(x_{ij}\right)$ be a random $m \times n$ matrix with
density $\psi$ of Section 2.5 , and let $p: {\Bbb R}^{mn}_+ \longrightarrow {\Bbb R}_+$
be the function defined in  Section 2.3. Then for some constant $\delta = \delta(\alpha, s_0)>0$ and 
some absolute constant $\kappa>0$, we have
$$\PP \left\{X: \quad p(X) > N^{\delta \ln N} \right\} \leq \kappa \left(2^{- m}+ 2^{-n}\right).$$

Therefore, the algorithm of Section~1.1 achieves $N^{O(\ln N)}$ complexity on these 
classes of margins.
\endproclaim

\proclaim{(3.3) Theorem} Fix $\alpha \geq 1$, $0 < \beta \leq 1$, and $\rho \geq 1$.
Suppose that $m \leq \rho n$, $n \leq \rho m$ and let $(R,C)$ be lower $\beta$ and
strongly upper $\alpha$-smooth margins.
Let $X=\left(x_{ij}\right)$ be a random $m \times n$ matrix with density
 $\psi$ of Section 2.5 and let $p: {\Bbb R}^{mn}_+ \longrightarrow {\Bbb R}_+$
be the function defined in  Section 2.3. Then for some constant $\delta = 
\delta(\rho, \alpha, \beta)>0$ and 
some absolute constant $\kappa>0$, we have
$$\PP \left\{X: \quad p(X) >   N^{\delta \ln N}\right\} \leq \kappa \left(2^{-m} +2^{-n}\right).$$

Therefore, the algorithm of Section~1.1 achieves $N^{O(\ln N)}$ complexity on these 
classes of margins.
\endproclaim

We remark that in 
Theorem 3.2 and Theorem 3.3 above, we can replace base 2 by any base $M>1$, 
fixed in advance.

\example{(3.4) Example:  symmetric margins} While conditions for $r_+$, $c_+$, $r_-$, and $c_-$ 
 are straightforward to verify, to check the upper bounds for $x_{ij}^{\ast}$ 
one may have to solve the optimization problem~(1.2.1) first. There are, however, some interesting
cases where an upper bound on $x_{ij}^{\ast}$ can be inferred from symmetry considerations.

Note that if two row sums $r_{i_1}$ and $r_{i_2}$ are equal then the transportation polytope 
$\PPP(R,C)$ is invariant under the transformation which swaps the $i_1$-st and $i_2$-nd
rows of a matrix $X \in \PPP(R,C)$. Since the function $g$ in the
optimization problem~(1.2.1) also remains 
invariant if the rows are swapped and is strictly concave, we must have $x_{i_1j}^{\ast}=x_{i_2j}^{\ast}$
for all $j$. Similarly, if $c_{j_1}=c_{j_2}$ we must have $x_{ij_1}^{\ast}=x_{ij_2}^{\ast}$ for 
all $i$. In particular, if all row sums are equal, we must have $x_{ij}^{\ast}=c_j/m$. Similarly, if 
all column sums are equal, we must have $x_{ij}^{\ast}=r_i/n$. 

More generally, 
one can show (see the proof of Theorem 3.5 in Section 6) that the largest 
entry $x_{ij}^{\ast}$ of $X^{\ast}$ necessarily lies at the intersection of 
 the row with the largest row sum $r_+$ 
and the column with the largest column sum $c_+$. Therefore, if 
$k$ of the row sums 
$r_i$ are equal to $r_+$ we must have $x_{ij}^{\ast} \leq c_+/k$. Similarly, if 
$k$ of the column sums are equal to $c_+$, we must have $x_{ij}^{\ast} \leq r_+/k$.
 
Here are some examples of classes margins where our algorithm provably achieves an $N^{O(\ln N)}$
complexity. 

$\bullet$ The class of margins for which at least a constant fraction of the row sums $r_i$ are equal to $r_+$:
$$\#\bigl\{i:\  r_i =r_+ \bigr\} =\Omega(m)$$
while  $m, n$, the row,  and the column sums 
differ by a factor, fixed in advance: $m/n=O(1)$, $n/m=O(1)$, $r_+/r_-=O(1)$, $c_+/c_-=O(1)$.
Indeed, in this case we have 
$$\max_{ij} x_{ij}^{\ast} = O(c_+/m)=O(N/mn)$$
and quasi-polynomiality follows by Theorem 3.3.

$\bullet$ The class of margins for which at least a constant fraction of the row sums $r_i$ are equal to $r_+$, while the column sums 
exceed the number of rows by at most a factor, fixed in advance, $c_+=O(m)$, and 
$m$ and $n$ are not too disparate: $m \leq 2^n$ and $n \leq 2^m$. Indeed, in this 
case 
$$\max_{ij} x_{ij}^{\ast} =O(c_+/m) =O(1)$$
and quasi-polynomiality follows by Theorem 3.2. 

$\bullet$ The classes of margins defined as above, but with rows and columns swapped.
\endexample

\medskip
For a different source of examples, we prove that if both ratios $r_+/r_-$ and 
$c_+/c_-$ are not too large, the margins 
are strongly upper smooth. To do this, we use the 
following general result about the typical table $X^{\ast}$, to be proved in Section~6:

\proclaim{(3.5) Theorem}

 Let $X^{\ast}=\left(x_{ij}^{\ast}\right)$ be the typical table.
\roster
\item We have 
$$x_{ij}^{\ast} \ \geq \ {r_- c_- \over r_+ m} \quad \text{and} \quad 
x_{ij}^{\ast} \ \geq \  {c_- r_- \over c_+ n} \quad \text{for all} \quad i,j.$$
\item If $r_-c_+ + r_- c_- + mr_- > r_+c_+$ then
$$x_{ij}^{\ast} \leq {c_+ \left(r_- c_- + m r_+\right) \over m \left( r_-c_+ + r_- c_- + mr_-  - r_+ c_+\right)}
\quad \text{for all} \quad i,j.$$
Similarly, if  
$c_- r_+ + c_- r_- + n c_- > r_+ c_+$ then 
$$x_{ij}^{\ast} \leq {r_+ \left(c_- r_- + n c_+\right) \over n \left( c_-r_+ + c_- r_- + nc_-  - c_+ r_+\right)}
\quad \text{for all} \quad i,j.$$
\endroster
\endproclaim

\example{(3.6) Example: golden ratio margins} Fix
$$1 \ \leq \ \beta \ <  \ {1 +\sqrt{5} \over 2} \approx 1.618$$
and a number $\rho \geq 1$.
Consider the class of margins $(R,C)$ such that $m \leq \rho n$, $n \leq \rho m$, and
$$r_+/r_-, \ c_+/c_- \ \leq \ \beta.$$ 
We claim that our algorithm has an $N^{O(\ln N)}$ complexity on this class
of margins.

To see this, let 
$$\beta_1=r_+/r_- \quad \text{and} \quad \beta_2=c_+/c_-.$$
If $\beta_1 \leq \beta_2$ then 
$$r_- c_+ + r_- c_-  - r_+ c_+ =\left(1+\beta_2 -\beta_1 \beta_2 \right) r_- c_-  \geq 
\left(1+\beta_2 -\beta_2^2 \right)r_- c_- \ \geq \ \epsilon r_-  c_-$$
for some $\epsilon =\epsilon(\beta)>0$ and hence by Part (2) of Theorem 3.5 we have 
$$x_{ij}^{\ast} \leq {c_+ \over m} \left({1 \over \epsilon} +\beta \right).$$
Similarly, if $\beta_2 \leq \beta_1$ then 
$$c_- r_+ + c_-r_- c_+ r_+  =\left(1 + \beta_1 -\beta_1 \beta_2 \right) r_-c_- \geq
\left(1+\beta_1 -\beta_1^2 \right) r_- c_- \ \geq \ \epsilon r_- c_-$$
for some $\epsilon=\epsilon(\beta)>0$
and hence
$$x_{ij}^{\ast} \leq {r_+ \over n} \left({1 \over \epsilon} +\beta \right).$$
In either case, $(R,C)$ are strongly upper $\alpha$-smooth for some 
$\alpha=\alpha(\beta)$ and Theorem 3.3 implies that our algorithm 
has a quasi-polynomial complexity on such margins. More generally, the algorithm 
is quasi-polynomial on the class of margins for which $\beta_1=r_+/r_-$ and 
$\beta_2=c_+/c_-$ are bounded above by a constant fixed in advance 
and $\beta_1 \beta_2 \leq \max\{ \beta_1, \beta_2\} +1-\epsilon$ where $\epsilon >0$
is fixed in advance.
\endexample

\example{(3.7) Example: linear margins} Fix $\beta \geq 1$ and $\epsilon>0$ 
such that $\epsilon \beta <1$
and consider the class of margins $(R, C)$ for which
$$r_+/r_- \leq \beta \quad \text{and} \quad c_+ \leq \epsilon m.$$
Part (2) of Theorem 3.5 implies that the margins $(R,C)$ are strongly upper $\alpha$-smooth 
for some $\alpha=\alpha(\beta, \epsilon)$ and therefore 
quasi-polynomiality of the algorithm is guaranteed by Theorem 3.2.
\endexample

The remainder of this paper is devoted to the proofs of
Theorems 3.2, 3.3, and 3.5. While the proof of Theorem 3.5
is relatively straightforward, our proofs of Theorem 3.2 and especially
Theorem 3.3 require some preparation. A 
general plan of the proofs of Theorems~3.2 and~3.3 is given in Section~5.

\head 4. Asymptotic estimates  \endhead

The following result proved in \cite{B07b} provides an asymptotic estimate for the number $\#(R, C)$ of 
contingency tables. It explains the role played by the optimization problem~(1.2.1). It will also introduces
ingredients needed in the statement and proof of 
Theorem~5.3 given below.

\proclaim{(4.1) Theorem} Let $\PPP(R,C)$ be the transportation
polytope of non-negative matrices with row sums $R$ and column sums $C$
and let $X^{\ast}=\left(x_{ij}^{\ast}\right)$ be the typical table, 
that is, the matrix $X^{\ast} \in \PPP(R,C)$ maximizing
$$g(X)=\sum_{ij} \Bigl( (x_{ij}+1)\ln (x_{ij}+1) -x_{ij} \ln x_{ij} \Bigr)$$
on $\PPP(R,C)$.
Let 
$$\rho(R, C)=\exp\left\{g(X^{\ast})\right\} = \max \Sb X =\left(x_{ij}\right) \\ X \in \PPP(R,C) \endSb
\prod_{ij} {\left(x_{ij}+1\right)^{x_{ij}+1} \over x_{ij}^{x_{ij}}}.$$
Then 
$$\rho(R, C) \  \geq\  \#(R, C) \ \geq \ N^{-\gamma(m+n)} \rho(R,C),$$
where $\gamma>0$ is an absolute constant.

Another representation of $\rho(R, C)$ is 
$$\rho(R, C) =\min \Sb 0 < x_1, \ldots, x_m <1 \\ 0 < y_1, \ldots, y_n <1 \endSb
\left(\prod_{i=1}^m x_i^{-r_i} \right) \left( \prod_{j=1}^n y_j^{-c_j} \right) 
\left(\prod_{ij} {1 \over 1-x_i y_j}\right).$$
A point $x_1, \ldots, x_m; y_1, \ldots, y_n$ minimizing the 
above product exists and is
unique up to scaling $x_i \longmapsto x_i\tau$, 
$y_j \longmapsto y_j\tau^{-1}$. It is related to $X^{\ast}$ by 
$$x_{ij}^{\ast}={x_i y_j \over 1-x_i y_j} \quad \text{for all} \quad 
i,j.$$
\endproclaim

We need the notion of the weighted enumeration of tables, as 
introduced in \cite{Ba08} and \cite{B07a}.

\subhead (4.2) Weighted enumeration of tables \endsubhead 
Fix margins $R$ and $C$ and  a non-negative $m \times n$ 
matrix $W$. Define
$$T(R,C; W)=\sum \Sb D=\left(d_{ij}\right) \endSb \prod_{ij} w_{ij}^{d_{ij}},$$
where the sum is taken over all $m \times n$ non-negative integer matrices $D$ with the 
row sums $R$ and  column sums $C$ and we agree that $w_{ij}^0=1$. Therefore,
$$\#(R,C)=T(R,C; {\bold1}),$$
where ${\bold 1}$ is the matrix of all 1's. 

The estimates of Theorem 4.1 extend to  weighted enumeration. We state only the 
part we are going to use. The following result is proved in \cite{B07b}.
\proclaim{(4.3) Theorem} Let 
$$\rho(R,C; W) =\inf \Sb x_1, \ldots, x_m > 0 \\ y_1, \ldots, y_n > 0 \\ 
w_{ij} x_i y_j <1 \ \text{for all} \  i,j \endSb \left(\prod_{i=1}^m x_i^{-r_i} \right) \left( \prod_{j=1}^n y_j^{-c_j} \right) \left(\prod_{ij} {1 \over 1-w_{ij} x_i y_j}\right).$$
Then
$$\rho(R,C; W)  \ \geq \ T(R,C; W) \ \geq \  N^{-\gamma(m+n)} \rho(R,C; W),$$
where $\gamma>0$ is an absolute constant.
\endproclaim
In fact, we will only use the upper bound of Theorem 4.3, 
which is actually straightforward to prove since 
$\prod_{ij} \left(1-w_{ij} x_i y_j \right)^{-1}$ is the generating function 
for the family $T(R,C; W)$.

\head 5. The plan of the proofs of Theorems 3.2 and 3.3 \endhead 

To prove Theorems 3.2 and 3.3 we need to understand the behavior of the function 
$$p(X) ={N^N \over N!} \per B(X),$$ 
that is, to estimate
values of permanents of doubly stochastic matrices. The following straightforward corollary of  
results of \cite{Fa81}, \cite{Eg81}, \cite{Br73},  and \cite{So03} shows that 
the permanent of an $N \times N$ doubly stochastic matrix lies close to 
$N!/N^N$ provided the entries of the matrix are not too large. We recall the definition of 
the Gamma function
$$\Gamma(t)=\int_0^{+\infty} x^{t-1} e^{-x} \ dx \quad \text{for} \quad t >0.$$

\proclaim{(5.1) Lemma} Let $B=\left(b_{ij}\right)$ be an $N \times N$ doubly stochastic 
matrix and let 
$$z_i=\max_{j=1, \ldots, N} b_{ij} \quad \text{for} \quad i=1, \ldots, N.$$
Suppose that
$$\sum_{i=1}^N z_i \leq \tau  \quad \text{for some} \quad \tau \geq 1.$$
Then 
$$\split {N! \over N^N} \ \leq\  \per B \ \leq \ &\left({ \tau \over N} \right)^N \Gamma^{\tau} \left(1+{N \over \tau} \right) \\ \leq \ & {N! \over N^N} \left(2 \pi N \right)^{\tau/2} e^{\tau^2/12 N} . 
\endsplit$$
\endproclaim
We delay the proof of Lemma 5.1 until Section 7.

We will apply Lemma 5.1 when $\tau = O(\ln N)$, in which case the ratio between the upper 
and lower bounds becomes $N^{O(\ln N)}$. In addition, 
we apply the lemma to the matrix $B(X)$, the doubly stochastic scaling of 
the random matrix $A(X)$ constructed in Theorem 2.1, see also Section 2.3. However,
to use this lemma, we need to bound 
the entries of $B(X)$. To do that, we will need to be able to bound the entries of the matrix 
$Y$ obtained from scaling $X$ to have
row sums $R$ and column sums $C$. To this end, we prove the following result
in Section~8, 
which might be of independent interest. 

\proclaim{(5.2) Theorem}  Let $R=\left(r_1, \ldots, r_m \right)$ and 
$C=\left(c_1, \ldots, c_n \right)$ be positive vectors such that 
$$\sum_{i=1}^m r_i = \sum_{j=1}^n c_j =N.$$
Let $X=\left(x_{ij}\right)$ be an $m \times n$ positive matrix 
and let $Y=\left(y_{ij}\right)$ be the scaling of $X$ to have row sums 
$R$ and column sums $C$, where 
$$y_{ij}=\lambda_i \mu_j x_{ij} \quad \text{for all} \quad i,j$$
and some positive $\lambda_1, \ldots, \lambda_m; \mu_1, \ldots, \mu_n$.

Then, for every $1 \leq p \leq m$ and $1 \leq q \leq n$ we have
$$\split  \ln y_{p q} \leq  \ln {r_p c_q \over N} + & \ln x_{pq} \\ +&\ln \left({1 \over N^2} \sum_{ij} r_i c_j x_{ij} \right) \\
 - &{1 \over N} \sum_{j=1}^n c_j \ln x_{pj} 
-{1 \over N} \sum_{i=1}^m r_i \ln x_{iq}.  \endsplit $$
\endproclaim 

Now suppose that $(R,C)$ are upper $\alpha$-smooth margins, that is, $r_i/N \leq \alpha/m$ and 
$c_j/N \leq \alpha/n$ for some $\alpha \geq 1$, fixed in advance. To give an idea of the
remainder of the argument and the role of the hypotheses, 
suppose further that $x_{ij}$ are sampled independently at random from the uniform 
distribution on $[0,1]$. Then Theorem 5.2 and the law of large numbers clearly
imply that as $m$ and $n$ grow, with overwhelming probability we have
$$y_{ij} \leq \kappa {r_i c_j \over N} x_{ij} \quad \text{for all} \quad i,j$$ 
and some absolute constant $\kappa >1$. 
If we construct the doubly stochastic matrix $B(X)$ as in Section 2.3, then with overwhelming 
probability for the entries $b_{ij}$ we will have 
$$b_{ij} \leq {\kappa \over N} \quad \text{for all} \quad i,j.$$

However, in the situation of our proof, the matrix $X=\left(x_{ij}\right)$ is actually
sampled from the distribution with density 
$\psi$ of Section 2.5. Thus to perform a similar analysis, we need to show that 
the entries of a random matrix $X$ are uniformly small. For that, we have 
to assume that the margins $(R,C)$ are {\it strongly} upper $\alpha$-smooth (in fact, one can 
show that merely the condition of upper smoothness is not enough). Specifically,
in Section 9, we prove 
the following result:

\proclaim{(5.3) Theorem} Let 
$$S \subset \Bigl\{ (i,j): \ i=1, \ldots, m; \quad j=1, \ldots, n \Bigr\}$$
be a set of indices, and let $X=\left(x_{ij}\right)$ be a random $m \times n$ matrix 
with density $\psi =\psi_{R,C}$ of Section 2.5. Suppose that the typical table 
$X^{\ast}=\left(x_{ij}^{\ast}\right)$ satisfies 
$$x_{ij}^{\ast} \leq \lambda \quad \text{for all} \quad i,j$$
and some $\lambda>0$.

Then for all $t >0$ we have 
$$\PP \left\{ \sum \Sb (i,j) \in S \endSb x_{ij} \geq t \right\} \leq
 \exp\left\{ -{t \over 2 \lambda+ 2} \right\} 4^{\#S} N^{\gamma(m+n)},$$
 where $\gamma>0$ is the absolute constant of Theorem 4.1.
\endproclaim

In Section 10 
we complete the proof of Theorem 3.2. Theorem 3.3 requires some more 
work and its proof is given in Section 12, after some technical estimates in Section 11.

\head 6. Proof of Theorem 3.5 \endhead

First, we observe that the typical table $X^{\ast}=\left(x_{ij}^{\ast}\right)$ 
is strictly positive, that is, it lies in the interior of the transportation polytope $\PPP(R,C)$.
Indeed, suppose that $x_{11}^{\ast}=0$, for example. Choose indices $p$ and $q$ such that
$x_{1q}^{\ast} >0$ and $x_{p1}^{\ast}>0$. Then necessarily $x_{pq}^{\ast} < r_p, c_q$
and we can consider a perturbation $X(\epsilon) \in \PPP(R,C)$ of $X^{\ast}$ defined 
for sufficiently small $\epsilon>0$ by 
$$x_{ij}=\cases x_{ij}^{\ast} +\epsilon &\text{if} \quad  i=1 \quad \text{and} \quad j=1 \\
                             x_{ij}^{\ast} -\epsilon &\text{if} \quad i=p, j=1 \quad \text{or} \quad i=1,j=q \\
                             x_{ij}^{\ast}+ \epsilon &\text{if} \quad i=p \quad \text{and} \quad j=q \\
                             x_{ij}^{\ast} &\text{if} \quad i \ne p \quad \text{and} \quad j \ne q. \endcases$$
Since the value of 
$${\partial  \over \partial x_{ij} } g(X) = \ln \left({x_{ij}+1\over x_{ij}}\right)$$
is equal to $+\infty$ at $x_{ij}=0$ (we consider the right derivative in this case)
 and finite if $x_{ij}>0$, we conclude that 
for a sufficiently small $\epsilon>0$, the matrix $X(\epsilon)$ attains a larger value of 
$g(X)$, which is a contradiction. We conclude that all the entries of the 
typical table $X^{\ast}$ are 
strictly positive.

Since $X^{\ast}$ lies in the interior of the transportation polytope
$\PPP(R,C)$, the Lagrange multiplier condition implies that
$$\ln \left({x_{ij}^{\ast}+1 \over x_{ij}^{\ast}} \right)=\lambda_i + \mu_j \quad \text{for all} \quad i,j 
\tag6.1$$
and some $\lambda_1, \ldots, \lambda_m$ and $\mu_1, \ldots, \mu_n$.
It follows that if $x_{i_1j}^{\ast} \geq x_{i_2 j}^{\ast}$ for some row indices $i_1, i_2$ and 
{\it some} column index $j$ then $\lambda_{i_1} \leq \lambda_{i_2}$ and hence 
$x_{i_1 j}^{\ast} \geq x_{i_2 j}^{\ast}$ for the same row indices $i_1$ and $i_2$ and {\it all} 
column indices $j$. 

We prove Part (1) first. Let us choose a row $i_0$ with the largest row sum $r_+$. 
Without loss of generality, we assume that $i_0 =1$. Hence 
$$x_{1j}^{\ast} \geq x_{ij}^{\ast} \quad \text{for} \quad j=1, \ldots, n.$$
Therefore, 
$$x_{1j}^{\ast}\  \geq \ {c_j \over m} \ \geq \ {c_- \over m} \quad \text{for} \quad j=1, \ldots, n.$$
Let us compare the entries in the first row and in the $i$-th row. From 
(6.1) we have 
$$\ln \left({x_{1j}^{\ast} +1 \over x_{1j}^{\ast} }\right) - \ln\left({x_{ij}^{\ast} +1 \over x_{ij}^{\ast}}\right)=
\lambda_1 - \lambda_i \quad \text{for} \quad j=1, \ldots, n. \tag6.2$$ 
Since 
$$\sum_{j=1}^n x_{1j}^{\ast} = r_+ \quad \text{and} \quad 
\sum_{j=1}^n x_{ij}^{\ast} \geq r_-,$$
there exists $j$ such that 
$${x_{ij}^{\ast} \over x_{1j}^{\ast}} \geq {r_- \over r_+}.$$
We apply (6.2) with that index $j$.
We have 
$$\lambda_1 - \lambda_i =
\ln {\left(x_{1j}^{\ast} +1\right)  x_{ij}^{\ast} \over \left(x_{ij}^{\ast}+1\right) x_{1j}^{\ast}}. \tag6.3$$
Now, the minimum value of 
$${(a+1) b \over (b+1) a} \quad \text{where}  \quad a \geq b \geq \tau a \quad 
\text{and} \quad a \geq \sigma$$
is attained at $a=\sigma$ and $b=\tau \sigma$ and equal to 
$${\tau \sigma + \tau \over \tau \sigma +1}.$$
In our case (6.3), 
$$ a=x_{1j}^{\ast}, \quad b=x_{ij}^{\ast}, \quad
\sigma = {c_- \over m}, \quad \tau = {r_- \over r_+}, \quad \text{and} \quad
{\tau \sigma + \tau \over \tau \sigma +1} ={r_- c_- + m r_- \over r_- c_- + m r_+}.$$
Hence 
$$\lambda_1 -\lambda_i \ \geq \ \ln {r_- c_- + m r_-  \over r_- c_- +m r_+ }.$$

Therefore, for every $j$,
$$\split \ln\left({x_{ij}^{\ast} +1 \over x_{ij}^{\ast}}\right) =&\ln \left({x_{1j}^{\ast}+1 \over x_{1j}^{\ast}} \right)
-\left(\lambda_1 -\lambda_i \right) \\ \leq &\ln \left({x_{1j}^{\ast}+1 \over x_{1j}^{\ast}} \right) -
\ln  {r_- c_- + m r_-  \over r_- c_- + mr_+ } \\ \leq & \ln {c_- +m \over c_-} - \ln {r_- c_- + m r_-  \over r_- c_- +m r_+ }. \endsplit $$
Hence
$${x_{ij}^{\ast} +1 \over x_{ij}^{\ast}} \ \leq\  {r_- c_- +r_+ m \over r_- c_-} \quad 
\text{for} \quad j=1, \ldots, n$$
and
$$x_{ij}^{\ast} \geq {r_- c_- \over r_+ m}$$
as desired. The second inequality in Part (1) is proved similarly.

To prove Part (2), we use an approach similar to that for Part (1), as well
as its inequality. 
Let $i_0$ be the row such that $r_{i_0}=r_-$. 
Without loss of generality,
we assume that $i_0=1$ and hence 
$$x_{ij}^{\ast} \geq x_{1j}^{\ast} \quad \text{for} \quad j=1, \ldots, n.$$
Thus we have 
$$x_{1j}^{\ast} \ \leq \ {c_j \over m} \ \leq \  {c_+ \over m} \quad \text{for} \quad j=1, \ldots, n.$$
Next, we compare the entries of the $i$-th row of $X^{\ast}$ and the entries of the first row 
using (6.2).

Since
$$\sum_{j=1}^n x_{ij}^{\ast} \leq r_+ \quad \text{and} \quad \sum_{j=1}^n x_{1j}^{\ast} = r_-$$
there is $j$ such that 
$${x_{ij}^{\ast} \over x_{1j}^{\ast}}  \leq {r_+ \over r_-}.$$

We apply (6.3) with that index $j$.
The maximum value of 
$${(a+1) b \over (b+1) a} \quad \text{where} \quad a \leq b \leq \tau a \quad \text{and} \quad 
a \geq \sigma$$
is attained at $a=\sigma, b=\tau \sigma$ and is equal to 
$${\tau \sigma + \tau \over \tau \sigma +1}.$$
In our case of (6.3),
$$a=x_{1j}^{\ast}, \quad b =x_{ij}^{\ast}, \quad
\tau = {r_+ \over r_-}, \quad \sigma ={ r_- c_- \over r_+ m}, \quad \text{and} \quad 
{\tau \sigma + \tau \over \tau \sigma +1} ={r_- c_- + m r_+ \over r_- c_- + m r_-} $$
where the expression for $\sigma$ follows by Part (1).
Hence
$$\lambda_1-\lambda_i \leq \ln {r_- c_- + m r_+ \over r_-c_- + m r_-}$$
and for all $j$ we have 
$$\split \ln\left({x_{ij}^{\ast} +1 \over x_{ij}^{\ast}}\right) =&\ln \left({x_{1j}^{\ast}+1 \over x_{1j}^{\ast}} \right)
-\left(\lambda_1 -\lambda_i \right) \\ \geq &\ln \left({x_{1j}^{\ast}+1 \over x_{1j}^{\ast}} \right) -
\ln  {r_- c_- +  m r_+ \over r_- c_- + m r_-} \\ \geq & \ln {c_+ +m \over c_+} - \ln {r_- c_- + m r_+ \over r_- c_- +m r_-}. \endsplit $$
Hence
$${x_{ij}^{\ast} +1 \over x_{ij}^{\ast}} \geq {\left(c_+ +m \right) \left(r_- c_- + m r_-\right) \over 
c_+ \left(r_- c_- + mr_+ \right)} \quad \text{for} \quad j=1, \ldots, n$$
and the proof follows.
{\hfill \hfill \hfill} \qed

\head 7. Proof of Lemma 5.1 \endhead 

We will use the following bounds for the permanent.

\subhead (7.1) The van der Waerden bound \endsubhead Let $B=\left(b_{ij}\right)$ be an 
$N \times N$ doubly stochastic matrix, that is, 
$$\sum_{j=1}^N b_{ij}=1 \quad \text{for} \quad i=1, \ldots, N \quad \text{and} 
\quad \sum_{i=1}^N b_{ij}=1 \quad \text{for} \quad j=1, \ldots, N$$
and 
$$b_{ij} \geq 0 \quad \text{for} \quad i,j=1, \ldots, N.$$
Then
$$ \per B \geq {N! \over N^N}.$$
This is the famous van der Waerden bound proved by Falikman \cite{Fa81} and Egorychev
\cite{Eg81}, see also Chapter 12 of \cite{LW01} and \cite{G06a}.

\subhead (7.2) The continuous version of the Bregman-Minc bound \endsubhead
Let $B=\left(b_{ij}\right)$ be an $N \times N$ matrix such that 
$$\sum_{j=1}^N b_{ij} \leq 1 \quad \text{for} \quad i=1, \ldots, N$$
and 
$$b_{ij} \geq 0 \quad i,j=1, \ldots, N.$$
Furthermore, let 
$$z_i =\max_{j=1, \ldots, N} b_{ij}>0 \quad \text{for} \quad i=1, \ldots, N.$$
Then 
$$\per B \leq \prod_{i=1}^N z_i \Gamma^{z_i} \left({1+z_i \over z_i}\right).$$
This bound was obtained by Soules \cite{So03}.

If $z_i=1/r_i$ for  integers $r_i$, the bound transforms into 
$$\per B \leq \prod_{i=1}^N {(r_i!)^{1/r_i} \over r_i},$$
which can be easily deduced from the Minc conjecture proved by  Bregman, see 
\cite{Br73}.

Now we are ready to prove Lemma 5.1.

\demo{Proof of Lemma 5.1}
The lower bound is the van der Waerden bound.

To prove the upper bound, define
$$f(\xi)=\xi\ln \Gamma \left({1+\xi \over \xi}\right)  +\ln \xi \quad \text{for} \quad 0< \xi \leq 1.$$
Then $f$ is a concave function and by the Bregman-Minc bound, we have 
$$\ln \per B \leq \sum_{i=1}^N f(z_i).$$
The function 
$$F(x) =\sum_{i=1}^N f(\xi_i) \quad \text{for} \quad x=\left(\xi_1, \ldots, \xi_N\right)$$ 
is concave on the simplex defined by the equation
$\xi_1 + \ldots + \xi_N = \tau$ and inequalities $\xi_i \geq 0$ for $i=1, \ldots, N$.
It is also symmetric under permutations of $\xi_1, \ldots, \xi_N$. 
Hence the maximum of $F$ is attained at 
$$\xi_1 = \ldots = \xi_N =\tau/N ,$$
and so 
$$\ln \per B \leq N f\left({\tau \over N} \right).$$
Thus
$$\per B \leq  \left({ \tau \over N} \right)^N \Gamma^{\tau} \left(1+{N \over \tau} \right)$$
and the rest follows by Stirling's formula.
{\hfill \hfill \hfill} \qed
\enddemo

\head 8. Proof of Theorem 5.2 \endhead 

We begin our proof by restating a theorem of Bregman \cite{Br73} in a slightly more general 
form.
\proclaim{(8.1) Theorem} Let $Y=\left(y_{ij}\right)$ be the positive $m \times n$ 
matrix that is the scaling of a positive $m \times n$ matrix $X=\left(x_{ij}\right)$ 
to have margins $(R,C)$.
Then
$$\sum_{ij} y_{ij} \left(\ln y_{ij} -\ln x_{ij}\right) \leq \sum_{ij} z_{ij}\left(\ln z_{ij}-\ln x_{ij}\right)$$
for every matrix $Z \in \PPP(R,C)$, where $\PPP(R,C)$ is the transportation polytope 
of $m \times n$ non-negative matrices with row sums $R$ and column sums $C$.
\endproclaim
\demo{Proof} The function 
$$f(Z)=\sum_{ij} z_{ij} \left(\ln z_{ij} - \ln x_{ij} \right)$$ 
is strictly convex on $\PPP(R,C)$ and hence attains its unique minimum 
$Y'=\left(y_{ij}'\right)$ on $\PPP(R, C)$. As in the proof of Theorem 3.5 (see Section 6), 
we can show that $Y'$ is strictly positive, that is, $Y'$ lies in the relative interior of $\PPP(R,C)$.
 Writing the Lagrange multiplier conditions, 
we obtain 
$$\ln y_{ij}' - \ln x_{ij}  =\xi_i + \eta_j$$ 
for some $\xi_1, \ldots, \xi_m$ and $\eta_1, \ldots, \eta_n$.
Letting $\lambda_i=e^{\xi_i}$ and $\mu_j =e^{\eta_j}$ we obtain 
$$y_{ij}' =\lambda_i \mu_j x_{ij} \quad \text{for all} \quad i,j,$$
so in fact $Y'=Y$
as desired.
{\hfill \hfill \hfill} \qed
\enddemo

Next, we prove a lemma that extends a result of Linial, Samorodnitsky, and  Wigderson
\cite{L+00}.
\proclaim{(8.2) Lemma} Let
$R=\left(r_1, \ldots, r_m \right)$ and  $C=\left(c_1, \ldots, c_n \right)$ be 
positive vectors
such that 
$$\sum_{i=1}^m r_i = \sum_{j=1}^n c_j =N.$$
Let $X=\left(x_{ij}\right)$ be an $m \times n$ positive matrix such that 
$$\sum_{ij} x_{ij}=N$$ and 
let $Y=\left(y_{ij}\right)$ be the scaling of $X$ to have row sums $R$ and column sums $C$. Then
$$\sum_{ij} r_i c_j \ln y_{ij} \geq \sum_{ij} r_i c_j \ln x_{ij}.$$ 
\endproclaim
\demo{Proof} Since $Y$ is the limit of the sequence of matrices obtained from $X$ by  
repeated alternate scaling of the rows to have row sums $r_1, \ldots, r_m$ and of the 
columns to have 
column sums $c_1, \ldots, c_n$, cf., for example, Chapter 6 of \cite{BR97},
it suffices to show that when the 
rows (columns) are scaled, the corresponding weighted sums of the logarithms of the 
entries of the matrix can only increase. 

To this end, let $X=\left(x_{ij}\right)$ be a positive $m \times n$ matrix with the 
row sums $\sigma_1, \ldots, \sigma_m$ such that 
$$\sum_{i=1}^m \sigma_i=N$$
and let 
$Y=\left(y_{ij}\right)$ be the matrix obtained from $Y$ by scaling the rows to have 
sums $r_1, \ldots, r_m$. Hence,
$$y_{ij}=r_i x_{ij}/\sigma_i \quad \text{for all} \quad i,j.$$
Thus
$$\sum_{ij} r_i c_j \left(\ln y_{ij}-\ln x_{ij}\right) =
\sum_{j=1}^n c_j \left(\sum_{i=1}^m \left(r_i \ln r_i - r_i \ln \sigma_i\right) \right) \geq 0,$$
since the maximum of the function
$$\sum_{i=1}^m r_i \ln \xi_i$$
on the simplex 
$$\left\{\sum_{i=1}^m \xi_i=N \quad \text{and} \quad \xi_i \geq 0 \quad 
\text{for} \quad i=1, \ldots, m \right\}$$
is attained at $\xi_i=r_i$.

The scaling of columns is treated similarly.
{\hfill \hfill \hfill} \qed
\enddemo

\demo{Proof of Theorem 5.2} 
Without loss of generality, we assume that $p=q=1$. 

Define an $m \times n$ matrix $U=\left(u_{ij}\right)$ by 
$$u_{ij}={r_i c_j x_{ij} \over  T} \quad \text{for} \quad T={1 \over N} \sum_{ij} r_i c_j x_{ij}.  \tag8.3$$
We note that  the scalings of $U$ and $X$ to margins $(R,C)$ coincide 
and that  
$$\sum_{ij} u_{ij}=N.$$

By Theorem 8.1, the matrix $Y$ minimizes
$$\sum_{ij} z_{ij} \left( \ln z_{ij}-\ln u_{ij}\right),$$
over the set $\PPP(R,C)$ of $m \times n$ non-negative matrices $Z$ with row sums $R$ and the 
column sums $C$.

For a real $t$, let us define the matrix $Y(t)=\left(y_{ij}(t)\right)$ by
$$y_{ij}(t)=\cases y_{ij}+t&\text{if \ } i=j=1 \\
 y_{ij}-{c_j \over N-c_1} t  &\text{if\ } i=1, j\ne 1\\
 y_{ij} -{r_i  \over N-r_1} t &\text{if\ } i\ne 1, j=1\\
 y_{ij}+{r_i c_j  \over (N-r_1)(N-c_1)} t &\text{if\ } i \ne 1, j \ne 1.
\endcases$$

Then $Y(0)=Y$ and $Y(t) \in \PPP(R,C)$ for all $t$ sufficiently close to 0.
Therefore, 
$${d \over dt} f\left(Y(t)\right) \bigm|_{t=0}=0,$$
where 
$$f(Z)=\sum_{ij} z_{ij} \left( \ln z_{ij} - \ln u_{ij} \right).$$
Therefore,
$$\split &\ln y_{11}-\ln u_{11}  +1\\
                  &\qquad   -{1 \over N-c_1} \sum_{j \ne 1} c_j \left( \ln y_{1j} -\ln u_{1j}  +1 \right) \\
                   &\qquad  -{1 \over N-r_1}\sum_{i \ne 1} r_i \left(\ln y_{i1}-\ln u_{i1} +1 \right) \\
                  & \qquad  +{1 \over (N-r_1)(N-c_1)}\sum_{i, j \ne 1} r_i c_j \left( \ln y_{ij}- \ln u_{ij} + 1 \right) \\
                     &=0.\endsplit
 $$
Rearranging the summands,
$$\split &{N^2 \over (N-r_1)(N-c_1)}
 \left(\ln y_{11}-\ln u_{11} \right) \\
                  &\qquad   - {N \over (N-r_1) (N-c_1)}  
                    \sum_{j =1}^n c_j \left( \ln y_{1j} -\ln u_{1j} \right) \\
                   &\qquad  -{N \over (N-r_1)(N-c_1)}\sum_{i =1}^m r_i \left(\ln y_{i1}-\ln u_{i1}\right) \\
                  & \qquad  +{1 \over (N-r_1)(N-c_1)}\sum_{ij} r_i c_j \left( \ln y_{ij}- \ln u_{ij} \right) \\
                     &=0.\endsplit
 $$
On the other hand, by Lemma 8.2,
$$\sum_{ij} r_i c_j \left(\ln y_{ij}-\ln u_{ij}\right) \geq 0,$$
so we must have
$$N^2 \left(\ln y_{11}-\ln u_{11}\right) 
-N \sum_{j =1}^n c_j \left(\ln y_{1j}-\ln u_{1j}\right) -N \sum_{i=1}^m r_i  \left( \ln y_{i1}-\ln u_{i1} \right) 
\leq 0. $$

In other words,
$$\ln y_{11} \leq \ln u_{11} + {1 \over N} \sum_{j=1}^n c_j \left( \ln y_{1j} - \ln u_{1j} \right) + 
{1 \over N} \sum_{i=1}^m r_i \left( \ln y_{i1} -\ln u_{i1} \right).$$
Since 
$$\sum_{j=1}^n y_{1j} = r_1,$$
we have 
$$\sum_{j=1}^n c_j \ln y_{1j}  \leq \sum_{j=1}^n c_j \ln \left({c_j r_1 \over N} \right),$$
cf. the proof of Lemma 8.2.
Similarly, since 
$$\sum_{i=1}^m y_{i1} = c_1,$$
 we have 
$$\sum_{i=1}^m r_i \ln y_{i1} \leq \sum_{i=1}^m r_i \ln \left({r_i c_1 \over N} \right).$$ 
Substituting (8.3) for $U$, we obtain 
$$\ln y_{11} \leq \ln x_{11} + \ln \left(r_1 c_1\right) -
 \ln T + {1 \over N} \sum_{j=1}^n c_j \ln{ T \over N x_{1j}}
+{1 \over N} \sum_{i=1}^m r_i \ln{T \over N x_{i1}},$$
and the proof follows.
{\hfill \hfill \hfill} \qed
\enddemo

\head 9. Proof of Theorem  5.3  \endhead

Fix margins $(R,C)$, let $\psi=\psi_{R,C}$ be the density of Section 2.5, and let 
$X=\left(x_{ij}\right)$ be the random matrix distributed in accordance with the density 
$\psi$. We will need a lemma that connects linear functionals of $X$ with the 
weighted sums $T(R,C;W)$ of Section 4.2.

\proclaim{(9.1) Lemma} 
Let $\lambda_{ij}<1$ be real numbers.
\roster 
\item Let $W=\left(w_{ij}\right)$ be the 
$m \times n$ matrix of weights given by
$$w_{ij}=\left(1-\lambda_{ij}\right)^{-1} \quad \text{for all} \quad i,j.$$
Then
$$\EE \exp\left\{ \sum_{ij} \lambda_{ij}  x_{ij} \right\}=
{T(R,C; W) \over \#(R,C)} \prod_{ij} w_{ij};$$
\item We have 
$$\EE \prod_{ij} x_{ij}^{-\lambda_{ij}} ={1 \over \#(R,C)} \sum \Sb D=\left(d_{ij}\right) \endSb
\prod_{ij} {\Gamma\left(d_{ij}-\lambda_{ij}+1\right) \over \Gamma\left(d_{ij}+1\right)},$$
where the sum is taken over all $m \times n$ non-negative integer matrices $D=\left(d_{ij}\right)$
with row sums $R$ and column sums $C$.
\endroster
\endproclaim
\demo{Proof} Let us prove Part (1). We have 
$$\split \EE \exp\left\{ \sum_{ij} \lambda_{ij} x_{ij} \right\} =
&{1 \over \#(R,C)} \int_{{\Bbb R}^{mn}_+} \exp\left\{ -\sum_{ij} \left(1-\lambda_{ij} \right) x_{ij}  \right\} \\
&\qquad \times \sum \Sb D=\left(d_{ij}\right) \endSb \prod_{ij} {x_{ij}^{d_{ij}} \over d_{ij}!} \ dX \\
=& {1 \over \#(R,C)} \int_{{\Bbb R}^{mn}_+} \exp\left\{ -\sum_{ij} x_{ij} \right\}  \\
&\qquad \times \sum \Sb D=\left(d_{ij}\right) \endSb  \prod_{ij} {w_{ij}^{d_{ij}}x_{ij}^{d_{ij}} \over d_{ij}!}
\prod_{ij} w_{ij} \ d X \\
= &{T(R,C; W) \over \#(R,C)} \prod_{ij} w_{ij},
\endsplit$$
as desired.

Since
$$\psi(X) \prod_{ij} x_{ij}^{-\lambda_{ij}} ={1 \over \#(R,C)} \sum \Sb D =\left(d_{ij}\right) \endSb
\prod_{ij} {x_{ij}^{d_{ij}-\lambda_{ij}} \over d_{ij}!} e^{-x_{ij}},$$
the proof of Part (2) follows.
{\hfill \hfill \hfill} \qed
\enddemo

To prove Theorem 5.3 we need only Part (1) of the lemma, while Part (2) will be used 
later in the proof of Theorem 3.3.

\demo{Proof of Theorem 5.3} We use the Laplace transform method, see, for example, 
Appendix A of \cite{AS92}. We have 
$$\split \PP \left\{ \sum \Sb (i,j) \in S \endSb x_{ij} \geq t  \right\} =
&\PP \left\{ \exp\left\{ {1 \over 2\lambda +2}\sum \Sb (i,j) \in S \endSb x_{ij} \right\}  \geq 
\exp\left\{{t \over 2\lambda+2} \right\}\right\} \\
\leq &\exp\left\{-{t \over 2\lambda+2} \right\}
 \EE \exp \left\{ {1 \over 2\lambda+2} \sum \Sb (i,j) \in S \endSb x_{ij} \right\},\endsplit $$
by the Markov inequality.

 By Part (1) of Lemma 9.1,
$$\EE \exp \left\{ {1 \over 2\lambda+2}\sum \Sb (i,j) \in S \endSb x_{ij} \right\} =
{T(R, C; W) \over \#(R, C)} \left({2\lambda + 2 \over 2\lambda +1}\right)^{\#S},$$
where 
$$w_{ij} =\cases (2\lambda+2)/(2\lambda+1) &\text{if\ } (i, j) \in S \\ 1 &\text{if \ } (i,j) \notin S. \endcases$$
Clearly,
$$\left({2\lambda + 2 \over 2\lambda +1}\right)^{\#S} \leq 2^{\#S}.$$
To bound the ratio of $T(R, C; W)$ and $\#(R, C)$, we use Theorems 4.1 and 4.3.

Let $0< x_1, \ldots, x_m; y_1, \ldots, y_n < 1$ be numbers such that
$$\rho(R,C)=\left(\prod_{i=1}^m {x_i}^{-r_i} \right) \left(\prod_{j=1}^n {y_j}^{-c_j} \right) 
\left( \prod_{ij} {1 \over 1-x_i y_j} \right).$$
For the typical table $X^{\ast}=\left(x_{ij}^{\ast}\right)$ we have 
$$x_{ij}^{\ast} = {x_i y_j \over 1-x_i y_j} \leq \lambda \quad \text{for all} \quad i,j.$$
Therefore, 
$$x_i y_j ={x_{ij}^{\ast} \over 1+x_{ij}^{\ast}}  \leq {\lambda \over \lambda+1} \quad 
\text{for all} \quad i,j$$
and 
$$w_{ij} x_i y_j <1 \quad \text{for all} \quad i,j.$$
Then we have
$$\rho(R,C;W) \leq \left(\prod_{i=1}^m {x_i}^{-r_i} \right) \left(\prod_{j=1}^n {y_j}^{-c_j} \right) 
\left( \prod_{ij} {1 \over 1-w_{ij} x_i y_j} \right)$$
and 
$${\rho(R,C;W) \over \rho(R,C)} \leq \prod \Sb (i,j) \in S \endSb {1-x_i y_j \over 1-w_{ij} x_i y_j}=
\prod \Sb (i,j) \in S \endSb {1 \over 1+(1-w_{ij}) x_{ij}^{\ast}}.$$
Now
$${1 \over 1 + (1-w_{ij}) x_{ij}^{\ast}} \leq {2\lambda +1  \over \lambda +1}
\leq 2 \quad \text{for all} \quad (i,j) \in S$$  and hence 
$${\rho(R, C; W) \over \rho(R, C)} \leq 2^{\#S}.$$
Since 
$$T(R, C; W) \ \leq \ \rho(R, C; W) \quad \text{and} \quad \#(R, C) \ \geq\ \rho(R, C) N^{-\gamma(m+n)},$$
 the proof follows.
 {\hfill \hfill \hfill} \qed
\enddemo

We will need the following corollary.
\proclaim{(9.2) Corollary} Suppose that $m \geq n$ and that the typical table
$X^{\ast}=\left(x_{ij}^{\ast}\right)$ satisfies
$$x_{ij}^{\ast} \leq \lambda \quad \text{for all} \quad i,j$$
and some $\lambda>0$.
Let $X=\left(x_{ij}\right)$ be a random $m \times n$ matrix distributed in accordance with the density 
$\psi_{R,C}$, and let
$$u_i=\max \Sb j=1, \ldots, n  \endSb x_{ij}.$$
Then for some $\tau=\tau(\lambda)>0$ we have 
$$\PP \left\{ \sum_{i=1}^m u_i \geq (\lambda+1) \tau m \ln N \right\} \leq 4^{-m}.$$
\endproclaim
\demo{Proof} We apply Theorem 5.3 to each of the $n^m$ of subsets $S$ having exactly one 
entry in each row.
{\hfill \hfill \hfill} \qed 
\enddemo 

We will also use an unconditional bound on the 
sum of all the entries of $X$.

\proclaim{(9.3) Lemma} We have
$$\PP \left\{ \sum_{ij} x_{ij}  \geq 2 (N +mn) \right\} \leq \left({3 \over 4} \right)^{N +mn}$$
\endproclaim
\demo{Proof} As in the proof of Theorem 5.3, we have
$$\split \PP \left\{ \sum_{ij} x_{ij} \geq 2 (N+mn) \right\} =
&\PP\left\{ \exp\left\{{1 \over 2} \sum_{ij} x_{ij} \right\} \geq \exp\left\{N+mn\right\} \right\}\\
\leq &\exp\{-(N +mn)\} \EE  \exp\left\{{1 \over 2} \sum_{ij} x_{ij} \right\}  \endsplit$$
by Markov's inequality.
 By Lemma 9.1,
$$\split \EE \exp\left\{ {1 \over 2} \sum_{ij} x_{ij} \right\} =&{T(R,C; W) \over \#(R,C)} \prod_{ij} w_{ij}
\quad \text{where} \\
&\qquad w_{ij}=2 \quad \text{for all} \quad i,j \\
=&2^{N+mn}\endsplit$$
and the proof follows.
{\hfill \hfill \hfill} \qed
\enddemo

\head 10. Proof of Theorem 3.2 \endhead 

We start with a technical result.

\proclaim{(10.1) Lemma} Let $(R,C)$ be upper $\alpha$-smooth margins, so 
$r_i/N \leq \alpha/m$ and $c_j/N \leq \alpha/n$ for all $i,j$.
Let $X=\left(x_{ij}\right)$ be a random $m \times n$ matrix with density $\psi_{R,C}$ of 
Section 2.5.
Then for any real $\tau$
$$\split &\PP \left\{ {1 \over N} \sum_{j=1}^n c_j \ln x_{ij} \leq -\tau  \right\} \leq 2^{n}
 \exp\left\{-{n \tau \over 2 \alpha}\right\} \quad \text{and} \\
&\PP \left\{ {1 \over N} \sum_{i =1}^m r_i \ln x_{ij} \leq -\tau \right\} \leq  
2^{m}\exp\left\{-{m \tau \over 2 \alpha}\right\}. \endsplit$$
\endproclaim  
\demo{Proof} Let us prove the first inequality.  As in the proof of Theorem 5.3, we use the Laplace transform method.
We have 
$$\split  \PP \left\{ {1 \over N} \sum_{j=1}^n c_j \ln x_{ij} \leq -\tau  \right\} = 
&\PP \left\{ -{n \over 2 \alpha N} \sum_{j=1}^n c_j \ln x_{ij} \geq  {n \tau \over 2 \alpha}  \right\} \\
\leq & \exp\left\{ -{n \tau \over 2 \alpha} \right\}
\EE \exp\left\{ -{n \over 2 \alpha N} \sum_{j=1}^n c_j \ln x_{ij}  \right\} \\
= & \exp\left\{ -{n \tau \over 2 \alpha} \right\}
\EE \prod_{j=1}^n x_{ij}^{-\lambda_j} \quad \text{where} \quad
\lambda_j={nc_j \over 2\alpha N}. \endsplit$$
Since 
$$\lambda_j \leq {1 \over 2},$$
by Part (2) of Lemma 9.1 we deduce that 
$$\EE \prod_{j=1}^n x_{ij}^{-\lambda_j} \leq \left(\Gamma \left({1 \over 2}\right)\right)^{n}
\leq 2^{n}$$
(we observe that every term in the sum of Lemma 9.1 does not exceed $\Gamma^n(1/2)$).
The proof of the second inequality is identical.
{\hfill \hfill \hfill} \qed
\enddemo

\demo{Proof of Theorem 3.2} 
Without loss of generality, we assume that $m \geq n$. We recall that function $p(X)$ 
is computed as follows. Given a positive $m \times n$ matrix 
$X=\left(x_{ij}\right)$, we 
compute the scaling $Y=\left(y_{ij}\right)$ of $X$ to have row sums $R$ and the column 
sums $C$.  Then we compute the $N \times N$ block matrix $B(X)$ 
consisting of $mn$ blocks of sizes
$r_i \times c_j$ with the entries in the $(i,j)$-th block equal to $y_{ij}/r_ic_j$.
Thus $B(X)$ is a doubly stochastic matrix and 
$$p(X)={N^N \over N!} \per B(X),$$
cf. Section 2.

We are going to use Theorem 5.2 to bound the entries of $Y$.

By Lemma 9.3,
$$\PP\left\{ \sum_{ij} x_{ij} < 2(N +mn) \right\}  \geq 1- \left({3 \over 4} \right)^{N+mn}.$$
Since $N \leq  s_0 mn$, $r_i/N \leq \alpha/m$,  and $c_j/N \leq \alpha/n$ we conclude that 
for some $\kappa_1=\kappa_1(\alpha, s_0)=2\alpha^2(s_0+1)$ we have 
$$\PP\left\{ {1 \over N^2} \sum_{ij} r_i c_j x_{ij} < \kappa_1 \right\}  
\geq 1- \left({3 \over 4} \right)^{N+mn}. $$
From Lemma 10.1, for a sufficiently large $\kappa_2 =\kappa_2(\alpha)$, we have 
$$\split &\PP \left\{ {1 \over N} \sum_{j=1}^n c_j \ln x_{pj} > -\kappa_2 \right\} \geq 1-4^{-n} \quad
\text{for all} \quad p=1, \ldots, m \quad 
\quad \text{and} \\
& \PP \left\{ {1 \over N} \sum_{i=1}^m r_i \ln x_{iq} > -\kappa_2 \right\} \geq 1- 4^{-m}
\quad \text{for} \quad q=1, \ldots, n. \endsplit$$

Therefore, by Theorem 5.2, we have for some $\kappa=\kappa(\alpha, s_0)$
$$\PP\left\{ y_{pq} \leq {r_p c_q \over N} \kappa x_{pq} \quad 
\text{for all} \quad p,q \right\} \geq 1-\left({3 \over 4}\right)^{N+nm}
-m4^{-n} -n 4^{-m}.$$

Now, $B$  consists of $mn$ blocks, the $(p,q)$-th block filled by the entries $y_{pq}/r_p c_q$.
Therefore the probability that for all $i,j=1, \ldots N$ we have 
$$b_{ij} \leq {\kappa \over N} x_{pq} \quad \text{provided} \quad (i,j) \quad \text{lies in the $(p,q)$-th
block of $B$} \tag10.2$$
is at least 
$$1-\left({3 \over 4}\right)^{N+nm}
-m4^{-n}-n4^{-m}.$$
We now bound $\per B(X)$ using Lemma 5.1 and Corollary 9.2.

Let 
$$\split &z_i=\max_{j=1, \ldots, N} b_{ij} \quad \text{for} \quad i=1, \ldots N \quad \text{and let}
\\ &u_p=\max_{q=1, \ldots, m} x_{pq}. \endsplit$$
Then, from (10.2) we have
$$\sum_{i=1}^N z_i \leq {\kappa \over N}  \sum_{p=1}^m r_p u_p \leq {\alpha \kappa \over m} 
\sum_{p=1}^m u_p.$$
By Corollary 9.2, for some $\tau_1=\tau_1(\alpha, s_0)$, we have 
$$\PP \left\{ \sum_{p=1}^m u_m \leq \tau_1  m \ln N \right\} \geq 1-4^{-m}.$$
Thus for some $\tau=\tau(\alpha, s_0)$ we have 
$$\PP\left\{ \sum_{i=1}^N z_i \leq \tau \ln N \right\} \geq 1 -\left({3 \over 4}\right)^{N+mn} 
-m 4^{-n} -n 4^{-m} -4^{-m}$$
and the proof follows by Lemma 5.1.
{\hfill \hfill \hfill} \qed
\enddemo

The rest of the paper deals with the proof of Theorem 3.3. This requires  sharpening of  the estimates of Lemma 10.1. Roughly, we need to prove that with 
overwhelming probability
$$\split &{1 \over N} \sum_{j=1}^n c_j \ln x_{ij} \geq -\tau + \ln s \quad \text{and} \\
&{1 \over N} \sum_{i=1}^m r_i \ln x_{ij} \geq -\tau + \ln s \endsplit$$
for some constant $\tau=\tau(\alpha, \beta)$, where $s=N/mn$ is the average entry of the table.

\head 11. An estimate of a sum
over tables \endhead 

To sharpen the estimates of Lemma 10.1 we need a more careful estimate of the sum
in Part (2) of Lemma 9.1. In this section, we prove the following technical result.

\proclaim{(11.1) Proposition} 
Suppose that $(R, C)$ are lower $\beta$-smooth and upper $\alpha$-smooth margins and that 
$$s=N/mn \geq 1.$$
Let $\lambda_1, \ldots, \lambda_m \leq 1/2$ be numbers and let $l=\lambda_1 +\ldots +\lambda_m$.
Then, for $k <n$ we have
$${1 \over \#(R,C)} \sum \Sb D=\left(d_{ij}\right) \endSb \prod \Sb 1 \leq i \leq m \\ 1 \leq j \leq k \endSb
{\Gamma(d_{ij}-\lambda_i +1) \over \Gamma(d_{ij}+1)} \leq \delta^{km} N^{\gamma(m+n)} s^{-kl},$$
where the sum is taken over all non-negative integer matrices $D$ with row sums $R$
and column sums $C$, 
 $\delta=\delta(\alpha, \beta)>0$  and $\gamma$ is the absolute constant of Theorem 4.1.
\endproclaim

We start with computing a simplified version of this sum in a closed form.
\definition{(11.2) Definition} Let us fix positive integers $c$ and $m$. The {\it integer simplex}
$\Upsilon(m, c)$ is the set of all non-negative integer vectors $a=\left(d_1, \ldots, d_m \right)$ 
such that $d_1 + \ldots + d_m =c$.

Clearly,
$$\# \Upsilon(m, c) ={m + c -1 \choose m-1}.$$
\enddefinition
A sum over $\Upsilon(m,c)$ similar to that of Proposition 11.1  can be computed 
in a closed form.
\proclaim{(11.3) Lemma} Let $\lambda_i <1$, $i=1, \ldots, m$, be numbers and let
$l=\lambda_1 +\ldots +\lambda_m$. Then
$${1 \over \# \Upsilon(m,c)} \sum \Sb d_1, \ldots, d_m \geq 0 \\ d_1 + \ldots +d_m=c \endSb
\prod_{i=1}^m {\Gamma\left(d_i-\lambda_i+ 1\right) \over \Gamma \left(d_i +1\right)} =
{\Gamma(c+m-l) \Gamma(m) \over \Gamma(c+m) \Gamma(m-l)}
 \prod_{i=1}^m \Gamma\left(1-\lambda_i \right).$$
\endproclaim
\demo{Proof} Let us define a function $h_c$ on the positive orthant ${\Bbb R}^m_+$ by 
the formula
$$h_c(x)={\left(m-1\right)! \over \left(m+c-1\right)!}
 \left(\sum_{i=1}^m \xi_i \right)^c \exp\left\{-\sum_{i=1}^m \xi_i \right\}
\quad \text{for} \quad x=\left(\xi_1, \ldots, \xi_m\right) \in {\Bbb R}^m_+.$$
Since
$$\left(\sum_{i=1}^m \xi_i \right)^c =\sum \Sb d_1, \ldots, d_m \geq 0 \\ d_1 + \ldots + d_m = c \endSb 
{c! \over d_1! \cdots d_m!} \xi_1^{d_1} \cdots \xi_m^{d_m},$$
We can rewrite
$$h_c(x)={m+c-1 \choose m-1}^{-1} 
\sum \Sb d_1, \ldots, d_m \geq 0 \\  d_1 + \ldots + d_m = c \endSb 
\prod_{i=1}^m {\xi_i^{d_i} \over d_i!} e^{-\xi_i}.$$
Therefore,
$${1 \over \# \Upsilon(m,c)} \sum \Sb d_1, \ldots, d_m \geq 0 \\ d_1 + \ldots +d_m=c \endSb
\prod_{i=1}^m {\Gamma\left(d_i-\lambda_i+ 1\right) \over \Gamma \left(d_i +1\right)}=
\int_{{\Bbb R}^m_+} h_c(x) \prod_{i=1}^m \xi_i^{-\lambda_i} \ dx.$$
Let $Q  \subset {\Bbb R}^m_+$ be the simplex 
$\xi_1 + \ldots + \xi_m=1$ with the Lebesgue measure $dx$ normalized to the probability measure.
Since the function 
$$\left(\sum_{i=1}^m \xi_i \right)^c \prod_{i=1}^m \xi_i^{-\lambda_i}$$ is positive homogeneous 
of degree $c-l$, we can write  
$$\int_{{\Bbb R}^m_+} h_c(x) \prod_{i=1}^m \xi_i^{-\lambda_i} \ dx=
{\Gamma(c+m-l) \over \Gamma(m)} \int_Q h_c(x) \prod_{i=1}^m \xi_i^{-\lambda_i} \ dx  \tag11.3.1
$$
On the other hand,
$$h_c(x)={\Gamma(m) \over \Gamma(c+m)} h_0(x) \quad \text{for} \quad  x \in Q. \tag11.3.2$$
Using (11.3.1) with $c=0$, we deduce that
$$\split \int_Q h_0(x) \prod_{i=1}^m \xi_i^{-\lambda_i} \ dx =&{\Gamma(m) \over \Gamma(m-l)}
\int_{{\Bbb R}^m_+} h_0(x) \prod_{i=1}^m \xi_i^{-\lambda_i} \ dx \\ =&
{\Gamma(m) \over \Gamma(m-l)}\prod_{i=1}^m \int_0^{+\infty} \xi_i^{-\lambda_i} e^{-\xi_i} \ d\xi_i \\ =
&{\Gamma(m) \over \Gamma(m-l)}\prod_{i=1}^m \Gamma\left(1-\lambda_i\right).\endsplit$$
Now, from (11.3.1) and (11.3.2), we have 
$$\int_{{\Bbb R}^m_+} h_c(x) \prod_{i=1}^m \xi_i^{-\lambda_i} \ dx  =
{\Gamma(c +m-l) \Gamma(m) \over \Gamma(c+m) \Gamma(m-l)} 
\prod_{i=1}^m \Gamma\left(1-\lambda_i\right),$$
as desired.
{\hfill \hfill \hfill}\qed
\enddemo
We need an estimate. 
\proclaim{(11.4) Corollary} Suppose that $\lambda_i < 1/2$ for $i=1, \ldots, m$ and $c \geq \beta m$ for some $\beta >0$. Then 
 $${1 \over \# \Upsilon(m,c)} \sum \Sb d_1, \ldots, d_m \geq 0 \\ d_1 + \ldots +d_m=c \endSb
\prod_{i=1}^m {\Gamma\left(d_i-\lambda_i+ 1\right) \over \Gamma \left(d_i +1\right)}\leq \left({m \over c}\right)^l \delta^m  $$
for some constant $\delta=\delta(\beta)>0$, where $l=\lambda_1 + \ldots + \lambda_m$.
\endproclaim
\demo{Proof} The proof follows from Lemma 11.3.
{\hfill \hfill \hfill} \qed 
\enddemo

Fix margins 
$R=\left(r_1, \ldots, r_m \right)$ and $C=\left(c_1, \ldots, c_n \right)$ and a number $k \leq n$.
Pick, uniformly at random, a 
contingency table $D=\left(d_{ij}\right)$ with margins $(R,C)$ and consider
its submatrix $Z$ consisting of the first $k$ columns. Hence $Z$ is an $m \times k$
non-negative integer matrix with the column sums $c_1, \ldots, c_k$. We interpret 
$Z$ as a point in the product
$$\Upsilon=\Upsilon\left(m, c_1\right) \times \cdots \times \Upsilon \left(m, c_k\right)$$
of integer simplices. This process induces 
a certain distribution on the set $\Upsilon$ of non-negative integer $m \times k$ matrices with
the column sums $c_1, \ldots, c_k$. 
We want to compare this distribution with the uniform distribution. 
Lemma~11.5 below says that the probability to get any particular matrix
$Z \in \Upsilon$ 
cannot exceed the 
uniform probability by much if the margins $(R,C)$ are smooth.

Once we fix the $m \times k$ submatrix $Z$ consisting of the first $k$ columns
 of a table with margins $(R,C)$,
the complementary  $m \times (n-k)$ table has row sums $R'=R-R(Z)$, where $R(Z)$ is the vector of row sums of $Z$, and 
column sums $\overline{C}=\left(c_{k+1}, \ldots, c_n \right)$, the truncation of $C$. 
Hence the probability of obtaining a particular $Z \in \Upsilon$ is 
$${\#(R', \overline{C})\over \#(R,C)},$$
where the ratio is declared to be 0 if $R'$ is not non-negative. 

We prove the following estimate.

\proclaim{(11.5) Lemma} Consider margins $(R,C)$ satisfying the 
constraints of Proposition~11.1.
Fix $k \leq n$
 and let $\Upsilon$ be the set of all $m \times k$ 
non-negative integer matrices with the column sums $c_1, \ldots, c_k$.

Let $\overline{C}=\left(c_{k+1}, \ldots, c_n \right)$, 
choose $Z \in \Upsilon$ and set
$R'=R-R(Z)$, where $R(Z)$ is the vector of the row sums of $Z$.
Then
$$ {\#(R', \overline{C}) \over \#(R,C)} \leq {\delta^{km} N^{\gamma(m+n)} \over \# \Upsilon}$$
for some constant $\delta=(\alpha, \beta) >0$, where $\gamma>0$ is an absolute constant
from Theorem~4.1.
\endproclaim

\demo{Proof} Let $\rho(R,C)$ be the quantity of Theorem 4.1. Here we agree that
$\rho(R', \overline{C})=0$ if $R'$ has negative components and that 
``$\max$'' and ``$\min$'' are replaced by ``$\sup$'' and ``$\inf$'' respectively if $R'$ is non-negative
but has 0 components.

Let $0< x_1, \ldots, x_m <1$ and $0< y_1, \ldots, y_n<1$ be an optimal point in Theorem
4.1, so
$$\rho(R, C)=  \prod_{i=1}^m x_i^{-r_i} \prod_{j=1}^n y_j^{-c_j} \prod \Sb 1 \leq i \leq m \\ 1 \leq j \leq n \endSb {1 \over 1- x_i y_j}.$$
Then 
$$\split \rho(R', \overline{C})  \leq & \prod_{i=1}^m x_i^{-r_i'} \prod_{j=k+1}^{n} y_j^{-c_j} 
\prod \Sb 1 \leq i \leq m \\ k+1 \leq j \leq n \endSb {1 \over 1-x_i y_j}  \\ \leq
&\prod_{i=1}^m x_i^{-r_i} \prod_{j=1}^n y_j^{-c_j} \prod  \Sb 1 \leq i \leq m \\ k+1 \leq j \leq n \endSb 
{1 \over 1-x_i y_j}
\endsplit$$
and hence
$${\rho(R', \overline{C}) \over \rho(R, C)} \leq \prod \Sb 1 \leq i \leq m \\ 1\leq j \leq k \endSb
 (1- x_i y_j).$$ 
Now, by Part (1) of Theorem 3.5, the typical table
$X^{\ast} =\left(x_{ij}^{\ast} \right)$ satisfies 
$$x_{ij}^{\ast} = {x_i y_j \over 1-x_i y_j} \geq \delta_1s \quad \text{for all} \quad i,j,$$ 
and for some $\delta_1=\delta_1(\alpha, \beta)$. This 
implies that 
$$1-x_i y_j ={1 \over 1+ x_{ij}^{\ast}} \leq {1 \over 1+\delta_1 s} \quad \text{for all} \quad i,j.$$
Summarizing,
$${\rho(R', \overline{C}) \over \rho(R,C)} \leq 
\left({1 \over 1+\delta_1 s} \right)^{km}.$$
Now,
$$\split \# \Upsilon= &\prod_{j=1}^k {c_j + m-1 \choose m-1} \leq
\prod_{j=1}^k  {c_j + m \choose m}\\  \leq & \prod_{j=1}^k \left({c_j + m \over c_j}\right)^{c_j} 
\left({c_j + m \over m} \right)^m. \endsplit$$
We have
$$\left({c_j + m \over c_j} \right)^{c_j} \leq e^m.$$ 
Furthermore, since $c_j \leq \alpha s m$, we have 
$$\left({c_j +m \over m} \right)^m \leq (1+\alpha s)^m$$ and 
$$\#\Upsilon \  {\rho(R', \overline{C}) \over \rho(R, C)} 
\leq e^{km} \left(1+\alpha s \over 1+\delta_1 s \right)^{km} \leq \delta^{km}.$$
Since by Theorem 4.1 we have 
$$\#(R,C) \ \geq \ N^{-\gamma(m+n)} \rho(R, C) \quad \text{and} \quad \#(R', \overline{C}) \ \leq \ 
\rho(R', \overline{C}),$$
the proof follows.
{\hfill \hfill \hfill} \qed
\enddemo

\demo{Proof of Proposition 11.1} Let $\Upsilon(m, c_j)$ be the integer simplex of non-negative integer 
vectors summing up to $c_j$ and let 
$$\Upsilon =\Upsilon(m, c_1) \times \cdots \times \Upsilon(m, c_k).$$
Using Lemma 11.5, we bound 
$$\split &{1 \over \#(R, C)} \sum \Sb D=\left(d_{ij}\right) \endSb \prod \Sb 1 \leq i \leq m \\ 1 \leq j \leq k \endSb 
 {\Gamma(d_{ij}-\lambda_i +1) \over \Gamma(d_{ij}+1)} \\ 
 =& \sum \Sb Z=\left(z_{ij} \right) \\ Z \in \Upsilon \endSb 
 {\#(R-R(Z), \ \overline{C}) \over \#(R, C)}  \prod  \Sb 1 \leq i \leq m \\ 1 \leq j \leq k \endSb 
 {\Gamma(z_{ij}-\lambda_i +1) \over \Gamma(z_{ij}+1)}
  \\ \leq
  &{\delta_1^{km} N^{\gamma(m+n)} \over \#\Upsilon}
  \sum \Sb Z=\left(z_{ij}\right) \\ Z \in \Upsilon \endSb  
  \prod \Sb 1 \leq i \leq m \\ 1 \leq j \leq k \endSb 
 {\Gamma(z_{ij}-\lambda_i +1) \over \Gamma(z_{ij}+1)} \endsplit$$  
for some $\delta_1=\delta(\alpha, \beta)$.
The sum
$${1 \over \# \Upsilon} \sum \Sb Z=\left(z_{ij}\right) \\ Z \in \Upsilon \endSb  
  \prod \Sb 1 \leq i \leq m \\ 1 \leq j \leq k \endSb 
 {\Gamma(z_{ij}-\lambda_i +1) \over \Gamma(z_{ij}+1)}$$
 is just the product of $k$ sums of the type
$${1 \over \Upsilon(m, c_j)} \sum  \Sb d_1, \ldots, d_m \geq 0 \\ d_1 +\ldots +d_m =c_j \endSb
\prod_{i=1}^m {\Gamma(d_i-\lambda_i +1) \over \Gamma(d_i+1)} \leq 
\left({m \over c_j}\right)^l \delta_2^m$$
by Corollary 11.4,
for some $\delta_2=\delta(\alpha, \beta)$.
The proof now follows.
{\hfill \hfill \hfill} \qed
\enddemo 

\head 12. Proof of Theorem 3.3 \endhead

Fix margins $(R,C)$ and let
 $X=\left(x_{ij}\right)$ be the $m \times n$ random matrix with density $\psi=\psi_{R,C}$ of Section 2.5. Define random variables
$$\split  &h_i={1 \over N} \sum_{j=1}^n c_j \ln x_{ij} \quad \text{for} \quad i=1, \ldots, m  \quad \text{and} \\
& v_j={1 \over N} \sum_{i=1}^m r_i \ln x_{ij} \quad \text{for} \quad j=1, \ldots, n. \endsplit$$
\proclaim{(12.1) Lemma} Let $(R,C)$ be lower $\beta$-smooth upper $\alpha$-smooth 
margins such that $s=N/mn \geq 1$. 

Choose a subset $J \subset \{1, \ldots, n\}$ of indices, $\#J=k$. 
Then for all $t>0$
we have 
$$ \PP\left\{ {1\over k} \sum \Sb j \in J \endSb v_j \leq -t + \ln s \right\} \leq \exp\left\{ -{tkm \over 2 \alpha} \right\} \delta^{km} N^{\gamma(m+n)},$$
Similarly, for  a subset $I \subset \{1, \ldots, m\}$ of indices, $\#I=k$,
we have 
$$ \PP\left\{ {1 \over k} \sum \Sb i \in I \endSb h_i \leq -t + \ln s  \right\} \leq 
\exp\left\{ -{tkn \over 2 \alpha} \right\} \delta^{kn} N^{\gamma(m+n)}.$$
for some number $\delta=\delta(\alpha, \beta)>0$  and the absolute constant $\gamma>0$ of 
Theorem 4.1.
\endproclaim
\demo{Proof} Without loss of generality, it suffices to prove only the first bound and only in the 
case of $J=\{ 1, \ldots, k \}$.

We use the Laplace transform method.
We have 
$$\split  \PP\left\{ {1 \over k} \sum_{j=1}^k v_j \leq -t + \ln s  \right\} =
&\PP\left\{ -{m \over 2 \alpha} \sum_{j=1}^k  v_j \geq {tkm  \over 2\alpha} - {km \ln s \over 2\alpha}  \right\}
\\ =&\PP \left\{ \exp\left\{ -{m \over 2 \alpha} \sum_{j=1}^k v_j \right\} \geq 
s^{-{km \over 2 \alpha}} \cdot \exp\left\{ {tk m \over 2 \alpha} \right\}  \right\} \\ 
\leq &s^{{km} \over {2\alpha}}\exp\left\{-{tkm \over 2 \alpha} \right\} 
\cdot  \EE\exp\left\{   -{m \over 2 \alpha} \sum_{j=1}^k  v_j  \right\}. \endsplit$$
 Let 
 $$\split &\lambda_i ={m r_i \over 2 \alpha N}  \leq {1 \over 2} \quad \text{for} \quad i=1, \ldots, m 
 \quad \text{and} \\
 &l=\lambda_1 + \ldots + \lambda_m = {m \over 2 \alpha}. \endsplit$$
 Using Part (2) of Lemma 9.1, we write
 $$ \EE\exp\left\{   -{m \over 2 \alpha} \sum_{j=1}^k  v_j  \right\} ={1 \over \#(R, C)} \sum \Sb D=\left(d_{ij}\right) \endSb \prod \Sb 1 \leq i \leq m \\ 1 \leq j \leq k \endSb 
 {\Gamma(d_{ij}-\lambda_i +1) \over \Gamma(d_{ij}+1)} ,  $$
 where the sum is taken over all contingency tables $D$ with margins $(R,C)$. 
 
The proof now follows by Proposition 11.1.
 {\hfill \hfill \hfill} \qed
\enddemo

We will use the following corollary.

\proclaim{(12.2) Corollary} 
Let $(R,C)$ be lower $\beta$-smooth upper $\alpha$-smooth margins such that $s=N/mn \geq 1$. 
Suppose further that $m \leq \rho n$ and $n \leq \rho m$ for some $\rho \geq 1$.Then 
for some $\tau=\tau(\alpha, \beta, \rho)>0$ we have 
$$\split &\PP\Bigl\{ \#\bigl\{i:\ h_i \leq -\tau  + \ln s \bigr\}  > \ln N \Bigr\} \leq 4^{-n} \quad \text{and}
\\ &\PP\Bigl\{ \#\bigl\{j:\ v_j \leq -\tau  + \ln s \bigr\}  > \ln N \Bigr\} \leq 4^{-m}. \endsplit$$
\endproclaim
\demo{Proof} We introduce random sets 
$$I=\bigl\{i: \ h_i \leq -\tau +\ln s \bigr\} \quad \text{and} \quad J=\bigl\{j: \ v_j \leq -\tau +\ln s \bigr\}$$
and note that 
$${1 \over \#I } \sum_{i \in I} h_i \leq -\tau +\ln s \quad \text{and} \quad {1 \over \#J } \sum_{j \in J} v_j \leq -\tau +\ln s.$$
The proof now follows from Lemma 12.1.
{\hfill \hfill \hfill} \qed
\enddemo

\demo{Proof of Theorem 3.3 } The proof is a modification of that of Theorem 3.2. We recall that
$$p(X)={N^N \over N!} \per B(X),$$
where $B(X)$ is the $N \times N$ doubly stochastic matrix constructed as follows: we scale
$m \times n$ matrix $X$ to the matrix $Y$ with row sums $R$ and column sums $C$
and let $b_{ij}=y_{pq}/r_p c_q$ provided the entry $(i,j)$ lies in the $(p,q)$-th block $B(X)$ of size
$r_p \times c_q$.
We are going to bound
the entries of  $Y$.
First, without loss of generality we assume that $s=N/mn \geq 1$ since the case of $s\leq 1$ is treated in
Theorem 3.2.

As in the proof of Theorem 3.2 we conclude that 
$$\PP\left\{ {1 \over N^2} \sum_{ij} r_i c_j x_{ij} < 2 \alpha^2 (s+1) \right\}  \geq 1- \left({3 \over 4} \right)^{N+mn}. \tag12.3$$

Let $$\split  &h_p={1 \over N} \sum_{j=1}^N c_j \ln x_{pj} \quad \text{for} \quad p=1, \ldots, m  \quad \text{and} \\
& v_q={1 \over N} \sum_{i=1}^m r_i \ln x_{iq} \quad \text{for} \quad q=1, \ldots, n. \endsplit$$
Choose $\tau>0$ as in Corollary 12.2. Set
$$P=\bigl\{p: \quad h_p \leq -\tau +\ln s \bigr\} \quad \text{and} \quad Q=\bigl\{q: \quad v_q \leq -\tau +\ln s \bigr\}.$$
Thus the probability that 
$$\#P \leq \ln N \quad \text{and} \quad \#Q \leq \ln N$$
is at least
$$1-4^{-m}-4^{-n}.$$
If $p \notin P$ and $q \notin Q$ and (12.3) holds then by Theorem 5.2, 
$$y_{pq} \leq \delta_1 {r_p c_q \over sN} x_{pq}$$
for some $\delta_1(\alpha, \beta)>0$.
If $ p \in P$  or $q \in Q$ then
$$y_{pq} \leq \min \{ r_p, c_q\}.$$

Consequently, for $b_{ij}$ with $(i,j)$ in the $p,q$-th block we have 
$$b_{ij} \leq   {\delta_1 \over sN} x_{pq} \quad \text{if} \quad p \notin P \quad \text{and}
\quad q \notin Q$$
and 
$$b_{ij} \leq \min\left\{{1 \over r_p}, \ {1 \over c_q} \right\} \quad \text{if} \quad p \in P
\quad \text{or} \quad q \in Q.$$

As in the proof of Theorem 3.2, we let 
$$\split &z_i=\max_{j=1, \ldots, N} b_{ij} \quad \text{for} \quad i=1, \ldots N \quad \text{and let}
\\ &u_p=\max_{q=1, \ldots, m} x_{pq}. \endsplit$$
We estimate that
$$z_i \leq {1 \over r_p}$$
if $i$ lies in the $p$-th row block with $p \in P$ and we estimate that 
$$ z_i \leq {\delta_1 \over sN} u_p +\max_{q \in Q}  {y_{pq} \over r_p c_q}, $$
if $i$ lies in the row block $p \notin P$.
Hence 
$$\sum_{i=1}^N z_i \leq \#P + {\delta_1 \over sN} \sum_{p=1}^m r_p u_p +
\sum_{p=1}^m \max_{q \in Q} {y_{pq} \over c_q}.$$
By Corollary 9.2,
$$\PP\left\{ \sum_{p=1}^m u_p \geq \tau_1 s m \ln N \right\} \leq 4^{-m}$$
for some $\tau_1=\tau_1(\alpha)$, and hence 
$$\PP \left\{{\delta_1 \over sN} \sum_{p=1}^m r_p u_p \leq \delta_2 \ln N \right\} \geq 1-4^{-m}.$$
for some $\delta_2=\delta_2(\alpha)$.
Finally,
$$\sum_{p=1}^m \max_{q \in Q} {y_{pq} \over c_q} \leq  \sum_{q \in Q} 
\sum_{p=1}^m {y_{pq} \over c_q} \leq \delta_3 \#Q$$
for some $\delta_3=\delta_3(\alpha, \beta)$.
Summarizing,
$$\PP \left\{ \sum_{i=1}^N z_i \leq \delta \ln N \right\} \geq 1- \left({3 \over 4} \right)^{N+mn} -
4^{-n} -2 \cdot 4^{-m}$$ 
for some $\delta=\delta(\alpha, \beta, \rho)>0$ and 
the proof is completed as in Theorem 3.2.
{\hfill \hfill \hfill}\qed
\enddemo

 \head Acknowledgments \endhead
 
 The authors are grateful to Jes\'us De Loera who computed some of the values 
 of $\#(R,C)$ for us using his {\tt LattE} code. The fourth author would like to thank 
 Radford Neal and Ofer Zeitouni for helpful
discussions.

 The research of the first author was partially supported by NSF Grant DMS 0400617.
The research of the third author was partially supported by ISF grant 039-7165.
The research of the first and third authors was also partially supported 
by a United States - Israel BSF grant 2006377.
The research of the fourth author was partially completed while he  was an NSF sponsored visitor
at the Institute for Pure and Applied Mathematics at UCLA, during April-June 2006. 
The fourth author was also partially supported by NSF grant 0601010 and an NSERC Postdoctoral
fellowship held at the Fields Institute, Toronto.

\enddocument
\end